\documentclass[lettersize,journal]{IEEEtran}
\usepackage{amsmath,amsfonts,amssymb}
\usepackage{algorithmic}
\usepackage{algorithm}
\usepackage{array}
\usepackage{textcomp}
\usepackage{stfloats}
\usepackage{url}
\usepackage{verbatim}
\usepackage{multirow}
\usepackage{graphicx}
\usepackage{cite}
\usepackage[capitalise]{cleveref}
\usepackage{epsfig}
\usepackage{textcomp}
\usepackage{xcolor}
\usepackage{cases}
\usepackage{subcaption}
\graphicspath{{Images/}}

\begin{document}

\title{Finite-Horizon LQR Control of the LWR Traffic Flow Model via Variable Speed Limits}

\author{Brian Block,~\IEEEmembership{Member,~IEEE,} and Stephanie Stockar,~\IEEEmembership{Member,~IEEE}%
\thanks{Received 4 June 2026. This work was supported by the National Science Foundation under NSF CAREER Award Grant 2042354. (\emph{Corresponding author: Brian Block.})}%
\thanks{Brian Block and Stephanie Stockar are with the Department of Mechanical and Aerospace Engineering, The Ohio State University, 201 W 19th Ave, Columbus, OH 43210, USA (email: block.168@osu.edu; stockar.1@osu.edu).}%
}

\markboth{IEEE Transactions on Control Systems Technology,~Vol.~\#, No.~\#, Month~2026}%
{Shell \MakeLowercase{\textit{et al.}}: A Sample Article Using IEEEtran.cls for IEEE Journals}

\IEEEpubid{0000--0000/00\$00.00~\copyright~2026 IEEE}

\maketitle

\begin{abstract}
This article presents a finite-horizon linear quadratic regulator (LQR) for the control of the first-order Lighthill-Whitham-Richards (LWR) traffic model with a triangular fundamental diagram. The in-domain control action is realized through variable speed limits implemented as a source term in the governing hyperbolic partial differential equation (PDE). Unlike previous infinite-horizon formulations, the proposed approach develops a finite-horizon LQR framework in which the state feedback function varies in both space and time. This requires the solution of a nonlinear Riccati PDE, which is done analytically using the parametric method of characteristics for both free-flow and congested traffic regimes. The resulting finite-horizon feedback law depends on space, time, and traffic regime, enabling regulation over prescribed time horizons. A sensitivity analysis is conducted to investigate the effects of the LQR parameters and terminal cost weight. Simulations on both circular and straight-road scenarios demonstrate the effectiveness of the proposed controller and show improved regulation performance relative to the corresponding infinite-horizon formulation.

\end{abstract}

\section{INTRODUCTION}
\IEEEPARstart{T}{raffic} congestion is a persistent challenge on highways, with significant implications for travel times and energy consumption \cite{schrank_urban_2021}. Addressing this issue demands effective controllers capable of alleviating traffic jams within specified time frames. Recent advancements in transportation technology, particularly in roadside infrastructure intelligence, offer promising avenues for better control of traffic dynamics.

Macroscopic models, which use aggregated variables instead of individual driver behavior, are often employed to capture the dynamics of large-scale highway traffic. In these representations, traffic is characterized by partial differential equations (PDEs) expressed in terms of density, velocity, and flow, drawing analogies with gas dynamics \cite{treiber_traffic_2013}. In contrast, microscopic car-following models focus on individual vehicles and their interactions. Despite these differences, macroscopic models have been shown to reproduce similar results to car-following models under various traffic conditions \cite{block_analysis_2025}. 

A common macroscopic model in traffic control is the Lighthill-Whitham-Richards (LWR) model, a scalar, hyperbolic PDE that balances accuracy and model complexity \cite{lighthill_kinematic_1955, richards_shock_1956, tumash_robust_2019, tumash_boundary_2022, yu_bilateral_2021}. The LWR model relies on an equilibrium flow-density relationship to determine vehicle velocity, typically represented using either a triangular fundamental diagram or Greenshield’s model \cite{treiber_traffic_2013}. For control and optimization of PDEs, the governing equation is often converted into a set of coupled ordinary differential equations (ODEs) \cite{christofides_robust_1998,atwell_reduced_2001,block_stabilization_2023}. In traffic control, this discretization to ODEs is popular in models such as METANET \cite{zegeye_integrated_2013,spiliopoulou_macroscopic_2014,pasquale_new_2018,d_frejo_macroscopic_2019}. In this work, instead the control is developed directly on the traffic flow PDE. 

Traffic flow can be controlled via boundary control, using ramp metering to modulate the flow entering a section of road. In \cite{pasquale_closed-loop_2018}, a proportional-integral controller for ramp metering was proposed, reducing traffic congestion with multiple on-ramps. Boundary control has also been employed to track desired time and space dependent density trajectories \cite{tumash_robust_2019,tumash_boundary_2022}, resulting in a reduction of the $\mathcal{L}_2$-norm of the error between desired and actual density. Furthermore, backstepping control has demonstrated its capability to decrease congestion and total travel time as well as track specific density profiles \cite{yu_bilateral_2021}. Reinforcement learning has also been utilized for boundary control in \cite{yu_reinforcement_2022}, where it was compared against both backstepping and PI control. In boundary control of traffic PDEs, there have been several works utilizing optimal control \cite{tumash_robust_2019,tumash_boundary_2022,li_optimal_2014,li_optimal_2014-1,liu_stochastic_2018}. However, boundary control has limitations, particularly in maintaining controllability under varying traffic regimes \cite{pisarski_nash_2016}. \IEEEpubidadjcol

Conversely, traffic flow can be regulated within the domain through variable speed limits (VSL). Adjusting speed limits along a road segment influences traffic flow and therefore mitigates congestion. Both on-road \cite{papageorgiou_effects_2008} and simulation \cite{papamichail_integrated_2008,allaby_variable_2007} studies have shown that implementing VSLs can avoid or delay congestion, enhance safety, and improve stability of traffic flow. In practice,  implementations of VSLs follow a rule based strategy where speed limits are selected based on average traffic speed, density measurements, and traffic volume \cite{allaby_variable_2007,khondaker_variable_2015,sadat_simulation-based_2017,block_lq-informed_2025}. VSL control has also been applied to discrete traffic models using model predictive control to minimize emissions and overall time spent in traffic \cite{liu_model_2017}. While VSL strategies have primarily been developed for discrete models \cite{papamichail_integrated_2008,liu_model_2017} or microscopic models \cite{allaby_variable_2007,khondaker_variable_2015}, recent efforts have extended its application to continuous PDE models for stabilizing desired density profiles \cite{karafyllis_feedback_2019-1, block_lq_2024}. However, the use of optimal in-domain control remains relatively unexplored, especially finite-time horizon control.

One approach for realizing optimal in-domain PDE control is through the use of a linear quadratic regulator (LQR). In \cite{aksikas_state_2007}, an LQ-feedback operator was developed by solving a matrix Riccati differential equation. Although the operator was spatially dependent, assuming an infinite time horizon resulted in a stationary solution. This method was applied to control the distributed jacket temperature of a fixed-bed reactor, modeled by a hyperbolic PDE, to achieve a desired chemical concentration \cite{aksikas_lq_2009} and later extended to a non-isothermal packed bed catalytic reactor, involving a coupled parabolic-hyperbolic PDE \cite{aksikas_optimal_2017}. In \cite{block_lq_2024}, an infinite horizon LQR was developed for regulating traffic density in a highway segment. In this approach, the optimal solution was only able to regulate the traffic within the free-flow domain, by defining an operator that was contingent on the Greenshield fundamental diagram. The approach was developed on the linearized LWR model and was shown to be accurate on the original nonlinear model as well. This same approach was used in \cite{block_lq-informed_2025} to create a rule-based controller that was tuned using real world data and applied to a highway case study. The control action in the LQR approach described in \cite{block_lq_2024} was unconstrained, so in \cite{block_constrained_2025} the LQR was combined with a barrier function in a quadratic program setup.

The primary drawback of the infinite-horizon LQR approach is its inability to guarantee that the control objectives are met within a given time window. While finite-horizon LQR problems have been successfully addressed for systems governed by ODEs \cite{suicmez_optimal_2014, tran_control_2017}, this has not yet been the case for systems described by PDEs. This limitation is particularly critical in traffic management, where controllers must alleviate traffic congestion within finite time frames. This article extends the previously developed infinite-horizon controller to a finite-horizon formulation and expands its applicability to both free-flow and congested traffic regimes. This is accomplished by deriving an analytical solution to the time- and space-dependent Riccati equation that arises from the finite-horizon optimal control problem. The resulting Riccati equation is an operator PDE \cite{aksikas_lq-optimal_2006,aksikas_lq_2009,aksikas_optimal_2017,block_lq_2024} that depends on both time and space, whereas the infinite-horizon formulation depends only on the spatial derivative. The analytical solution to the operator Riccati PDE is found via the method of characteristics for both the free-flow and congested regimes. The  resulting state feedback function is regime dependent, yielding a switching control law between free-flow and congested traffic conditions. The resulting controller is then evaluated in both a straight highway road scenario and a circular road scenario. Additionally, a sensitivity analysis is performed with respect to the LQR weighting matrices $Q$ and $R$ for the  infinite-horizon controller, and  the terminal state weight $S$  and horizon length $t_f$ for the finite-horizon controller.

This article has three main contributions:
\begin{enumerate}
    \item The previously developed VSL controller \cite{block_lq_2024} is extended to both free-flow and congested traffic regimes through the use of regime-dependent linear operators.  
    \item A finite-horizon optimal control formulation is developed for the LWR traffic model. The resulting  operator Riccati equation becomes a nonlinear PDE in both space and time, for which an analytical solution is obtained via the method of characteristics. The resulting state feedback function depends on time, space, and  traffic regime.
    \item The proposed controller is evaluated on representative circular and straight road traffic scenarios. In addition, a sensitivity    
    analysis is conducted to investigate the influence of the weighting matrices  $Q$ and $R$ on the infinite-horizon controller and the influence of the terminal state weight $S$ and horizon length $t_f$ on the finite-horizon controller. 
\end{enumerate}

The structure of the article is as follows. \cref{sec:model} describes the LWR traffic model and the its linear form as well as the introduction of variable speed limits. \cref{sec:control} explains the control formulation extension to both free flow and congested traffic together with the extension to finite-time control. \cref{sec:sensitivity} investigates the sensitivity of the control input to the control parameters on two simplified, yet representative scenarios. Finally, \cref{sec:conclusion} highlights the impact of this work and presents the conclusions and possible future work in this area.

\section{MODEL DESCRIPTION}\label{sec:model}
\subsection{Lighthill-Whitham-Richards Model}
\begin{figure}[thpb]
  \centering
  \includegraphics[width=1\columnwidth]{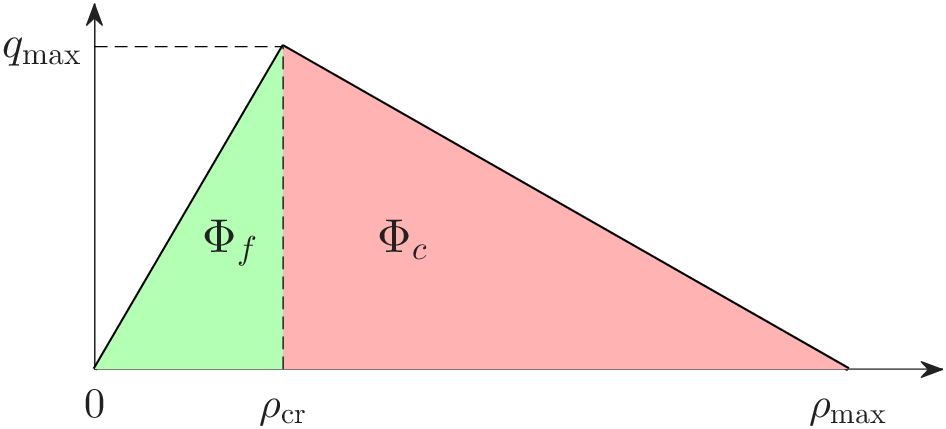}
  \caption{Triangular fundamental diagram.
  }
  \label{fig:fd}
\end{figure}

The LWR model \cite{lighthill_kinematic_1955, richards_shock_1956} is a first-order, hyperbolic conservation equation given by
\begin{equation} \label{eq:LWR}
    \frac{\partial\rho(z,t)}{\partial t} + \frac{\partial q(\rho(z,t))}{\partial z} = 0
\end{equation}
where $\rho(z,t)$ is the traffic density, representing the number of cars occupying a road length, $q(z,t)$ is the traffic flow rate, and $z$ and $t$ are the space and time variables, respectively. The closure of the problem is achieved through the adoption of a fundamental flow diagram that relates density with flow. The fundamental diagram used in this paper is the triangular model \cite{daganzo_requiem_1995}, shown in \cref{fig:fd} and described by
\begin{equation} \label{eq:fd}
    q(\rho(z,t)) = 
    \begin{cases} v_\mathrm{max}\rho(z,t), & \forall \rho(z,t)\in\Phi_f\\
                  w(\rho_\mathrm{max} - \rho(z,t)), & \forall \rho(z,t)\in\Phi_c        
    \end{cases}
\end{equation}
where $v_\mathrm{max}$ is the maximum velocity, $w$ is the congested speed, and $\rho_\mathrm{max}$ is the maximum density of the road. The free flow regime $\Phi_f$ is shown in green in \cref{fig:fd}, while the congested regime $\Phi_c$ is in red. The congested speed $w$ is given as
\begin{equation} \label{eq:w}
    w = \frac{v_\mathrm{max}\rho_\mathrm{cr}}{\rho_\mathrm{max}-\rho_\mathrm{cr}}
\end{equation}
where $\rho_\mathrm{cr}$ is the critical density where traffic moves from the free flow to the congested regime.

\subsection{Variable Speed Limits}
The variable speed limit actuation is included in the governing equations by modifying the fundamental diagram in \cref{eq:fd} as
\begin{equation} \label{eq:fdVSL}
    q(\rho(z,t),b(z,t)) = 
    \begin{cases} b(z,t)v_\mathrm{max}\rho(z,t), & \forall \rho\in\Phi_f\\
                  b(z,t)w(\rho_\mathrm{max} - \rho(z,t)), & \forall \rho\in\Phi_c        
    \end{cases}
\end{equation}
where $b(z,t)$ is the VSL ratio. A value of $b(z,t)>1$ increases the slope of both the free flow and congested regime in \cref{fig:fd}, while a value of $b(z,t)<1$ will decrease it. This, in turn affects the maximum possible flow $q_\mathrm{max}$ which occurs at the critical density. In this work, the location of $\rho_\mathrm{cr}$ is held constant. For simplicity, for the rest of the paper $(z,t)$ is omitted, except where it is needed for emphasis.

\subsection{Linearization}
The linearized LWR model is derived considering the nominal equilibrium density $\rho_0$ and nominal maximum speed ratio $b_0$ \cite{block_lq_2024} such that density, speed ratio, and flow can be written as
\begin{equation}\label{eq:linearizedVSLfd}
    \begin{aligned}
        \rho &= \rho_0 + \Delta\rho\\
        b &= b_0 + \Delta b\\
        q &=
        \begin{cases}
            v_\mathrm{max}((\rho_0 + \Delta\rho)(b_0 + \Delta b)), & \forall \rho\in\Phi_f\\
            w((\rho_\mathrm{max} - \rho_0 - \Delta\rho)(b_0 + \Delta b)), & \forall \rho\in\Phi_c
        \end{cases}
    \end{aligned}
\end{equation}
Then, the linearized LWR model can then be written as
\begin{equation}\label{eq:linearizedLWR}
    \frac{\partial\Delta\rho}{\partial t} + q_0\frac{\partial\Delta q}{\partial z} = 0
\end{equation}
In \eqref{eq:linearizedLWR}, $q$ is the generalized flow function. To account for the two distinct traffic regimes, namely free and congested flow, two separate models are developed with the corresponding flow functions for each regime:
\subsubsection{Free Flow}
When $\rho \leq \rho_\mathrm{cr}$, the linearized model is
\begin{equation} \label{eq:FFlinLWR}
    \frac{\partial\Delta\rho}{\partial t} + v_\mathrm{max}b_0\frac{\partial\Delta\rho}{\partial z} + v_\mathrm{max}\rho_0\frac{\partial\Delta b}{\partial z}= 0
\end{equation}
\subsubsection{Congested}
When $\rho > \rho_\mathrm{cr}$, the linearized model is
\begin{equation} \label{eq:ClinLWR}
    \frac{\partial\Delta\rho}{\partial t} - w b_0\frac{\partial\Delta\rho}{\partial z} + w(\rho_\mathrm{max}-\rho_0)\frac{\partial\Delta b}{\partial z} = 0
\end{equation}
In both cases, the variable $\frac{\partial\Delta b}{\partial z}$ represents the change in the VSL ratio over the length of road. If no VSL is applied, $\frac{\partial\Delta b}{\partial z} = 0$ and $b_0 = 1$, resulting in the standard linear LWR model, \cite{block_lq_2024}.

\section{CONTROL FORMULATION}\label{sec:control}
This section extends the previously developed infinite horizon LQR \cite{block_lq_2024} to manage both free flow and congested regimes by defining separate analytical state feedback functions for each regime. The LQR problem is then reformulated over a finite horizon, resulting in the state feedback function to become dependent on both space and time. The state feedback function is found via the method of parametric characteristics. This new optimal control formulation enables effective control within specified time frames which is paramount for traffic congestion control.

The theory considered \cite{aksikas_lq-optimal_2006,aksikas_lq_2009,block_lq_2024} considers a linear hyperbolic PDE of the form
\begin{equation} \label{eq:standardlinearPDE}
    \begin{aligned}
        \frac{\partial x}{\partial t}(z,t) &= V\frac{\partial x}{\partial z}(z,t) + Mx(z,t) + Bu(z,t) \\
        y(z,t) &= Cx(z,t)
    \end{aligned}
\end{equation}
The equivalent state-space form of \eqref{eq:standardlinearPDE} is given by the following differential equation on the Hilbert space $\mathcal{H}$
\begin{equation} \label{eq:SSHyper}
    \begin{aligned}
        \Dot{x} &= Ax +Bu\\
        y &= Cx
    \end{aligned}
\end{equation}
where $A$ is a linear operator defined by
\begin{equation} \label{eq:linoperator}
    A = V\cdot\frac{d.}{dz} + M\cdot I
\end{equation}
on the domain
\begin{equation} \label{eq:domain}
    D(A) = \{x\in\mathcal{H}:\frac{dx}{dz}\in\mathcal{H}\textrm{ and }x(0)=0\}
\end{equation}
If $V<0$, then the operator $A$ \eqref{eq:linoperator}-\eqref{eq:domain} generates an exponentially stable $C$-semigroup, ensuring that $(A,B)$ is exponentially stabilizable and $(C,A)$ is exponentially detectable \cite{aksikas_lq_2009}. Conversely, according to Remark 3 in \cite{aksikas_lq_2009} if $V>0$, then the exponential stability can be proven given that the boundary condition $x(0)=0$ in \eqref{eq:domain} becomes $x(1)=0$, so
\begin{equation} \label{eq:domain2}
    D(A) = \{x\in\mathcal{H}:\frac{dx}{dz}\in\mathcal{H}\textrm{ and }x(1)=0\}
\end{equation}
In order to define the operator $A$ \eqref{eq:linoperator}-\eqref{eq:domain2} for the LWR, \eqref{eq:linearizedLWR} can be written in the form of \eqref{eq:standardlinearPDE} where $V$, $M$, $B$, $C$, $x$, and $u$ are given as
\begin{equation} \label{eq:LWRLQRvariables}
    \begin{aligned}
        V &= \begin{cases} -v_\mathrm{max}b_0, & \forall \rho\in\Phi_f\\
                           wb_0, & \forall \rho\in\Phi_c        
             \end{cases}\\
        M &= 0\\
        B &= \begin{cases} -v_\mathrm{max}\rho_0, & \forall \rho\in\Phi_f\\
                           -w(\rho_\mathrm{max}-\rho_0), & \forall \rho\in\Phi_c        
             \end{cases}\\
        C &= I\\
        x &= \Delta\rho,\quad u = \frac{\partial\Delta b}{\partial z}
    \end{aligned}
\end{equation}

\subsection{Infinite Horizon LQR for LWR Model}
For the system given by \eqref{eq:SSHyper}-\eqref{eq:LWRLQRvariables}, consider the LQ optimal control problem on an infinite time interval, where for any initial condition $x_0\in\mathcal{H}$ the cost function is 
\begin{equation} \label{eq:costfunctINF}
    J(x_0,u_\mathrm{opt}) = \frac{1}{2}\int_{0}^{\infty} \langle x(t),Qx(t) \rangle + \langle u(t), Ru(t) \rangle \,dt
\end{equation}
where $u_\mathrm{opt}$ is the control input that minimizes \eqref{eq:costfunctINF}, $Q$ is a positive semi-definite matrix that penalizes deviations from the desired state, and $R$ is a positive definite matrix that penalizes control effort. The optimal control input $u_\mathrm{opt}$ is determined by finding the non-negative self-adjoint operator $P$ that solves the operator Riccati equation
\begin{equation} \label{eq:ORE}
    [A^*P + PA + Q - PBR^{-1}B^*P]x = 0
\end{equation}
for all $x\in D(A)$, where $P(D(A))\subset D(A^*)$. When $(A,B)$ is exponentially stabilizable and $(C,A)$ is exponentially detectable, \cref{eq:ORE} has a unique, non-negative solution $P$, ensuring that the cost functional \cref{eq:costfunctINF} is minimized \cite{aksikas_lq_2009}. When $V<0$ and $A$ is defined on the domain given by \eqref{eq:domain}, an exponentially stable semigroup is generated and these condition hold. Similarly, when $V>0$ and $A$ is defined on the domain given by \eqref{eq:domain2}, these conditions also hold. Additionally, as long as $V>0$ and $Q>0$ in \cref{eq:ORE}, then the solution to \cref{eq:ORE} is unique \cite{aksikas_optimal_2017}. Furthermore, $A$ is exponentially stable if either \emph{(i)} $V$ is diagonalizable with all identical eigenvalues or \emph{(ii)} the eigenvalues of the matrix $V$ are negative \cite{aksikas_lq-optimal_2006}. Given the solution $P$, the optimal control input is 
\begin{equation} \label{eq:optcont}
    u_\mathrm{opt}(z,t) = Kx(z,t)
\end{equation}
where for the infinite horizon problem the feedback $K=K_\mathrm{inf}$ is given as
\begin{equation}
    K_\mathrm{inf}=-R^{-1}B^*P(z)
\end{equation}
To determine the state feedback function $P$ the definition of $A$ is inserted into \eqref{eq:ORE} and the resulting equation is
\begin{equation}\label{eq:InfHorOde}
    V\frac{dP}{dz} = M^*P + PM + Q - PBR^{-1}B^*P
\end{equation}
In order to solve \eqref{eq:InfHorOde}, a boundary condition must be defined. When $V<0$ and the domain of $A$ is given by \eqref{eq:domain}, if $P(L)=0$, where $L$ is the length of the spatial domain, then $P(L)x(L)=0$ for all $x\in D(A)$ which implies $P(D(A))\subset D(A^*)$. On the other hand, if $V>0$ and the domain of $A$ is defined by \eqref{eq:domain2}, then $P(0)=0$ makes $P(0)x(0)=0$ for all $x\in D(A)$ which again implies $P(D(A))\subset D(A^*)$. So, the boundary condition for \eqref{eq:InfHorOde} can be given as
\begin{equation}\label{eq:InfHorBCs}
    \begin{cases}
        P(L)=0 & \mathrm{if} \quad V<0\\
        P(0)=0 & \mathrm{if} \quad V>0
    \end{cases}
\end{equation}
The analytical solution to \eqref{eq:InfHorOde} for $V<0$ will first be solved completely, after which the the solution will be stated for when $V>0$, which is a trivial extension of the same process. Using separation of variables, \eqref{eq:InfHorOde} becomes
\begin{equation}\label{eq:SolveInfHorP1}
    \int_{0}^P \frac{V}{Q-PBR^{-1}B^*P} dP = \int_{0}^z dz = z + c_1
\end{equation}
The left-hand side of \eqref{eq:SolveInfHorP1} can be solved by a combination of integration by partial fractions and substitution, taking into account that $Q$, $B$, and $R$ are scalars. After applying the aforementioned methods, \eqref{eq:SolveInfHorP1} becomes
\begin{equation}\label{eq:SolveInfHorP2}
    \frac{V\sqrt{R}}{2B\sqrt{Q}}\bigg(\ln{\bigg|P + \frac{\sqrt{QR}}{B}\bigg|} - \ln{\bigg|P - \frac{\sqrt{QR}}{B}\bigg|}\bigg) = z + c_1
\end{equation}
Then,
\begin{equation}\label{eq:SolveInfHorP3}
    \begin{aligned}
        \ln{\bigg|\frac{P + \frac{\sqrt{QR}}{B}}{P - \frac{\sqrt{QR}}{B}}\bigg|} &= \frac{2B\sqrt{Q}}{V\sqrt{R}}(z+c_1)\\
        \bigg|\frac{P + \frac{\sqrt{QR}}{B}}{P - \frac{\sqrt{QR}}{B}}\bigg| &= \exp{\bigg(\frac{2B\sqrt{Q}}{V\sqrt{R}}(z+c_1)\bigg)}\\
        \bigg|P + \frac{\sqrt{QR}}{B}\bigg| &= \bigg|P - \frac{\sqrt{QR}}{B}\bigg|\exp{\bigg(\frac{2B\sqrt{Q}}{V\sqrt{R}}(z+c_1)\bigg)}
    \end{aligned}
\end{equation}
The solution to \eqref{eq:SolveInfHorP3} can be either
\begin{equation}\label{eq:SolveInfHorP4_incorrect}
    P + \frac{\sqrt{QR}}{B} = \bigg(P - \frac{\sqrt{QR}}{B}\bigg)\exp{\bigg(\frac{2B\sqrt{Q}}{V\sqrt{R}}(z+c_1)\bigg)}
\end{equation}
or
\begin{equation}\label{eq:SolveInfHorP4}
    P + \frac{\sqrt{QR}}{B} = -\bigg(P - \frac{\sqrt{QR}}{B}\bigg)\exp{\bigg(\frac{2B\sqrt{Q}}{V\sqrt{R}}(z+c_1)\bigg)}
\end{equation}
Because we are solving for the solution when $V<0$ the boundary condition $P(L)=0$ only holds for \eqref{eq:SolveInfHorP4}. The final solution for the state feedback function $P$ using the condition $P(L)=0$ is then
\begin{equation}\label{eq:InfHorPSolution1}
    P(z) = \frac{\sqrt{QR}}{B}\cdot\frac{\exp{\bigg(\frac{2B\sqrt{Q}}{V\sqrt{R}}(z-L)\bigg)}-1}{\exp{\bigg(\frac{2B\sqrt{Q}}{V\sqrt{R}}(z-L)\bigg)}+1}
\end{equation}
If instead $V>0$ and $P(0)=0$ then the solution for the state feedback function becomes
\begin{equation}\label{eq:InfHorPSolution2}
    P(z) = \frac{\sqrt{QR}}{B}\cdot\frac{\exp{\bigg(\frac{2B\sqrt{Q}}{V\sqrt{R}}z\bigg)}-1}{\exp{\bigg(\frac{2B\sqrt{Q}}{V\sqrt{R}}z\bigg)}+1}
\end{equation}

\subsection{Finite Horizon LQR for LWR Model}
\begin{figure}[ht!]
    \centering
    \begin{subfigure}{0.49\columnwidth}
            \centering
            \includegraphics[width=\textwidth]{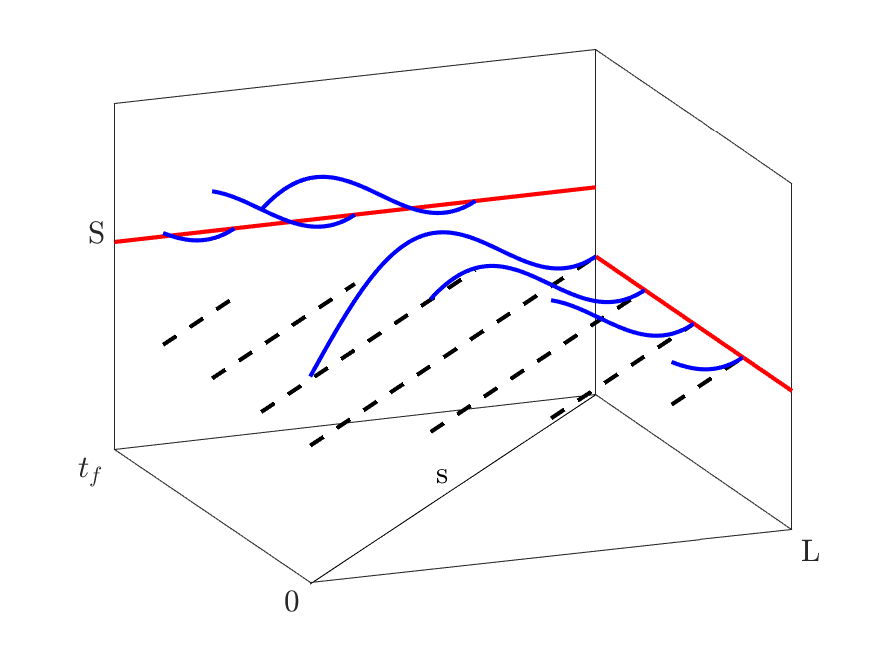}
            \caption{}
            \label{fig:characteristics}
    \end{subfigure}
    \begin{subfigure}{0.49\columnwidth}
            \centering
            \includegraphics[width=\textwidth]{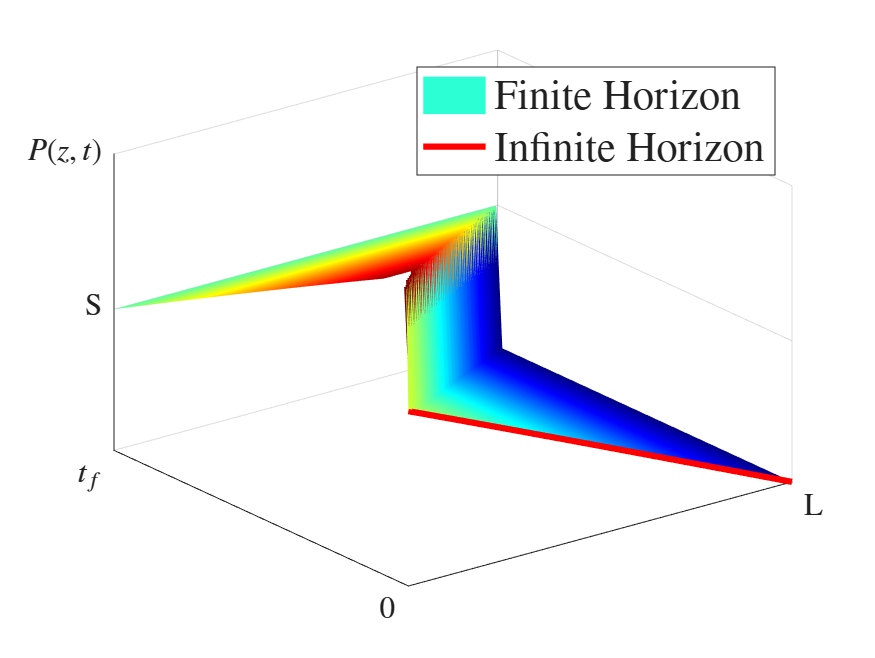}
            \caption{}
            \label{fig:statefeedback}
    \end{subfigure}
    \caption{Depiction of boundary and final conditions for method of characteristics (a) and comparison of infinite and finite horizon state feedback functions (b) in free-flow regime. }
    \label{fig:statefeedbackexamples}
\end{figure}

Now, the finite time horizon LQR problem is considered for the system \eqref{eq:SSHyper}-\eqref{eq:LWRLQRvariables}. The cost function for any initial condition $x_0\in\mathcal{H}$ is now written as
\begin{equation} \label{eq:costfunctFIN}
    J(x_0,u) = \frac{1}{2}\langle x(t_f),Sx(t_f) \rangle + \frac{1}{2}\int_{0}^{t_f} \langle x,Qx \rangle + \langle u, Ru \rangle \,dt
\end{equation}
where $S$ is a symmetric positive semi-definite matrix that places a weighted cost on the final state $x(t_f)$. Because of the final time constraint, the non-negative self-adjoint operator $P$ is now found by solving the differential Riccati equation
\begin{equation} \label{eq:DRE}
    \bigg[\frac{\partial P}{\partial t} + A^*P + PA + Q - PBR^{-1}B^*P\bigg]x=0
\end{equation}
Given the solution $P$ to \eqref{eq:DRE} and feedback in \eqref{eq:optcont} is now $K=K_\mathrm{fin}$ which is now a function of both space and time
\begin{equation}
    K_\mathrm{fin} = -R^{-1}B^*P(z,t)
\end{equation}
Once again, the definition of $A$ is inserted into the Riccati equation \eqref{eq:DRE}, resulting in
\begin{equation} \label{eq:FinHorPde}
    -\frac{\partial P}{\partial t} = -V \frac{\partial P}{\partial z} + M^*P +PM + Q - PBR^{-1}B^*P
\end{equation}
Unlike the infinite horizon case, this now becomes a PDE. In order to solve it, initial and/or final conditions must be placed on the value of $P$ in both the time and space domain. The boundary condition in space stays the same from the infinite horizon case but applied for all time, so
\begin{equation}\label{eq:FinHorBCs}
    \begin{cases}
        P(L,t)=0 & \mathrm{if} \quad V<0\\
        P(0,t)=0 & \mathrm{if} \quad V>0
    \end{cases}
\end{equation}
Then, to enforce finite time convergence, the condition placed on the value of $P$ at the final time $t_f$ is
\begin{equation}\label{eq:FinHorFC}
    P(z,t_f) = S
\end{equation}
Since \eqref{eq:FinHorPde} with conditions \eqref{eq:FinHorBCs}, \eqref{eq:FinHorFC} is a first-order PDE, it is solved analytically via the parametric form of the method of characteristics
\begin{equation}\label{eq:parametrics}
    z = z(s), \quad t = t(s), \quad P = P(s)
\end{equation}
where $s$ is the distance along the associated curve. Then, using \eqref{eq:FinHorPde} and \eqref{eq:parametrics}, the Lagrange-Charpit equations for the system can be written as
\begin{equation}\label{eq:lagrangecharpit}
    \frac{dz}{V} = \frac{dt}{-1} = \frac{dP}{Q - PBR^{-1}B^*P} = ds
\end{equation}
To find the solution, the following characteristic system of ODEs needs to be solved
\begin{subnumcases}{}
    \label{eq:characteristic_t}\frac{dt}{ds} = -1 & $t(0) = 0$\\
    \label{eq:characteristic_z}\frac{dz}{ds} = V  & $z(t_f(s)) = \xi$\\
    \label{eq:characteristic_P}\frac{dP}{ds} = Q - PBR^{-1}B^*P & $P(t_f(s)) = F(\xi)$
\end{subnumcases}
The solution to \eqref{eq:characteristic_t} and \eqref{eq:characteristic_z} is trivial, resulting in
\begin{equation}\label{eq:t_solution}
    t = -s
\end{equation}
\begin{equation}\label{eq:z_solution}
    z = V(t_f + s) + \xi
\end{equation}
The boundary conditions in the parameter $s$ for \eqref{eq:characteristic_z} and \eqref{eq:characteristic_P} change depending on where the characteristic crosses the boundary. This is illustrated for the free-flow regime in \cref{fig:characteristics} where the black dotted lines correspond to the characteristics, the blue lines correspond to the value of $P$ along the characteristics, and the red lines are the boundary values. Considering \eqref{eq:FinHorBCs} and \eqref{eq:FinHorFC}, if a characteristic curve hits the final time boundary before $z = L$, the value of of $F(\xi)$ is $S$. If instead the characteristic never reaches the final time boundary and instead reaches the condition at $z=L$ first, then $F(\xi) = 0$. Finally, \eqref{eq:characteristic_P} is solved in the same manner as the infinite-horizon LQR problem and can be written as
\begin{equation}\label{eq:P_solution}
    P = \frac{\sqrt{QR}}{B}\cdot\frac{1 - \exp{\bigg(-2B\sqrt{\frac{Q}{R}}(s+c_2)}\bigg)}{1 + \exp{\bigg(-2B\sqrt{\frac{Q}{R}}(s+c_2)}\bigg)}
\end{equation}
where
\begin{equation}\label{eq:c_2}
    \begin{aligned}
        c_2 &= -\frac{1}{2B}\sqrt{\frac{R}{Q}}\ln{\bigg|\frac{F(\xi)-\sqrt{QRB^{-2}}}{F(\xi)+\sqrt{QRB^{-2}}}\bigg|} + s(\xi)\\
        F(\xi) &= \begin{cases}
        S & \xi<L\\
        0 & \xi\geq L
        \end{cases}\\
        s(\xi) &= \begin{cases}
        -t_f & \xi<L\\
        \frac{1}{V}(L-\xi-Vt_f) & \xi\geq L
        \end{cases}\\
    \xi &= z - V(s+t_f)
    \end{aligned}
\end{equation}
An example of a finite horizon state feedback function $P(z,t)$ is shown in \cref{fig:statefeedback} where the surface depicts the finite horizon solution and the red line depicts the infinite horizon solution $P(z)$ for the same parameters in the free-flow regime. If the time horizon is long enough the finite time horizon solution starts out the same as the infinite horizon solution. 

\section{PARAMETER SENSITIVITY ANALYSIS}\label{sec:sensitivity}
\subsection{Impact of $Q$ and $R$ on Infinite-Horizon Optimal Solution}

\begin{figure}[ht!]
    \centering
    \begin{subfigure}{0.32\columnwidth}
            \centering
            \includegraphics[width=\textwidth]{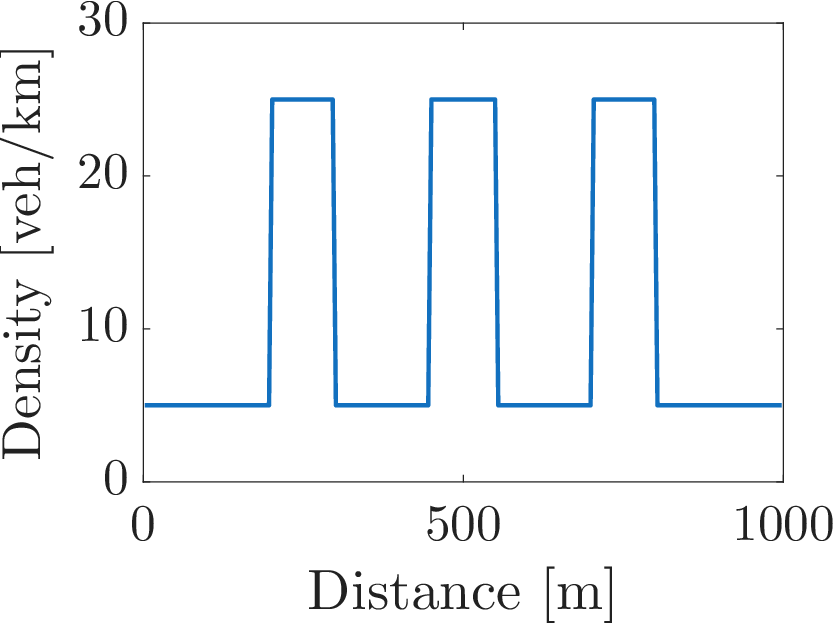}
            \caption{}
            \label{fig:CircularRoadIC}
    \end{subfigure}
    \begin{subfigure}{0.32\columnwidth}
            \centering
            \includegraphics[width=\textwidth]{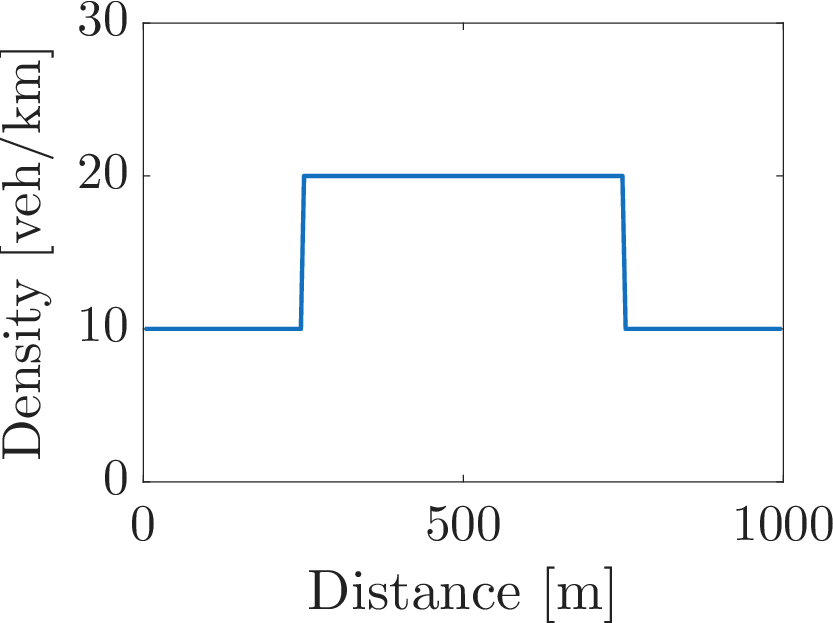}
            \caption{}
            \label{fig:StraightRoadIC}
    \end{subfigure}
    \begin{subfigure}{0.32\columnwidth}
            \centering
            \includegraphics[width=\textwidth]{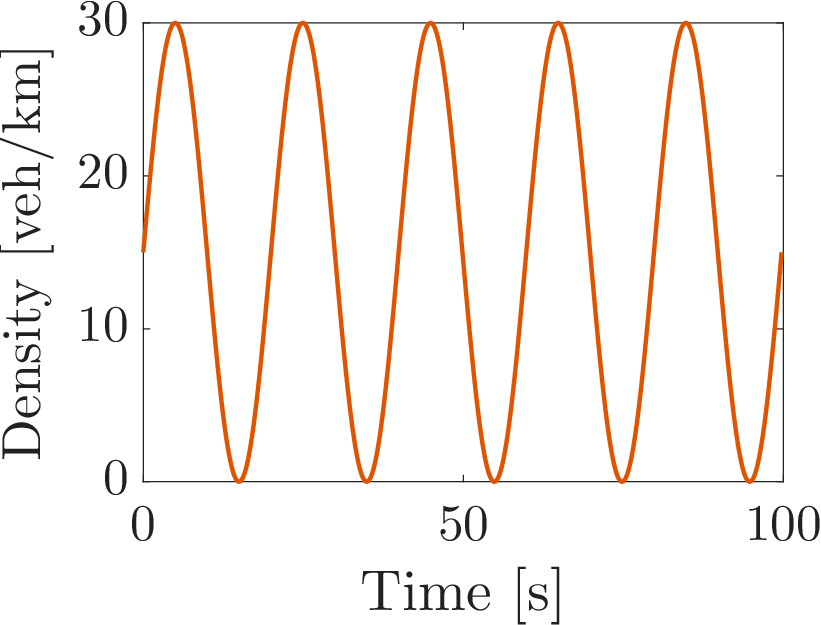}
            \caption{}
            \label{fig:StraightRoadBC}
    \end{subfigure}
    \caption{Initial density on circular road (a), initial density on straight road (b), and density inflow for straight road (c). }
    \label{fig:InitConditionsSimStudy}
\end{figure}

Because the controller is unconstrained, the effort to achieve the control goal varies greatly with different values of $Q$ and $R$. To analyze the impact of the LQR parameters, two different simulations are done: one on a circular road and one on a straight road. The parameters for both simulations are given in \cref{tab:params_QR_simulation_study}. The circular road is modeled with periodic boundary conditions such that
\begin{equation}\label{eq:circularroadperiodicBC}
    \rho(0,t) = \rho(L,t)
\end{equation}
and with an initial condition shown in \cref{fig:CircularRoadIC} and given by
\begin{equation}\label{eq:circularroadIC}
    \rho(z,0) = \begin{cases}
        5  & z\leq 200\\
        20 & 200<z\leq300\\
        5  & 300<z\leq450\\
        20 & 450<z\leq550\\
        5  & 550<z\leq700\\
        20 & 700<z\leq800\\
        5  & x> 800
    \end{cases}
\end{equation}
For the straight road, the traffic flow is allowed to flow freely out of the boundary at $z=L$ and the boundary condition at the inlet of the road is shown in \cref{fig:StraightRoadBC} given by
\begin{equation}\label{eq:straightroadBC}
    \rho(0,t) = 15 + 15\sin{\bigg(\frac{\pi t}{10}\bigg)}
\end{equation}
while the initial condition is shown in \cref{fig:StraightRoadIC} and given by
\begin{equation}\label{eq:straightroadIC}
    \rho(z,0) = \begin{cases}
        10 & z\leq 250\\
        20 & 250<z<750\\
        10 & x\geq 750
    \end{cases}
\end{equation}
In both case studies the traffic is mixed, meaning parts of the road are in free flow and parts are congested.

\begin{table}[!ht]
    \centering
    \begin{center}
    \caption{Model Parameters.}
    \label{tab:params_QR_simulation_study}
    \begin{tabular}{l|l|l }
    \hline \hline
    Parameter & Value & Unit\\  \hline 
    Maximum density $\rho_\mathrm{max}$ & 150 & [veh/km] \\
    Critical density $\rho_\mathrm{cr}$ & 15 & [veh/km] \\
    Maximum speed $V_\mathrm{max}$ & 72 & [km/h] \\
    Road length $L$ & 1000 & [m]\\
    Simulation time $T$ & 100 & [s]\\
     \hline \hline
    \end{tabular}
    \end{center}
\end{table}

\subsubsection{Circular Road Case Study}

\begin{figure*}[ht!]
    \centering
    \begin{subfigure}{0.32\textwidth}
            \centering
            \includegraphics[width=\textwidth]{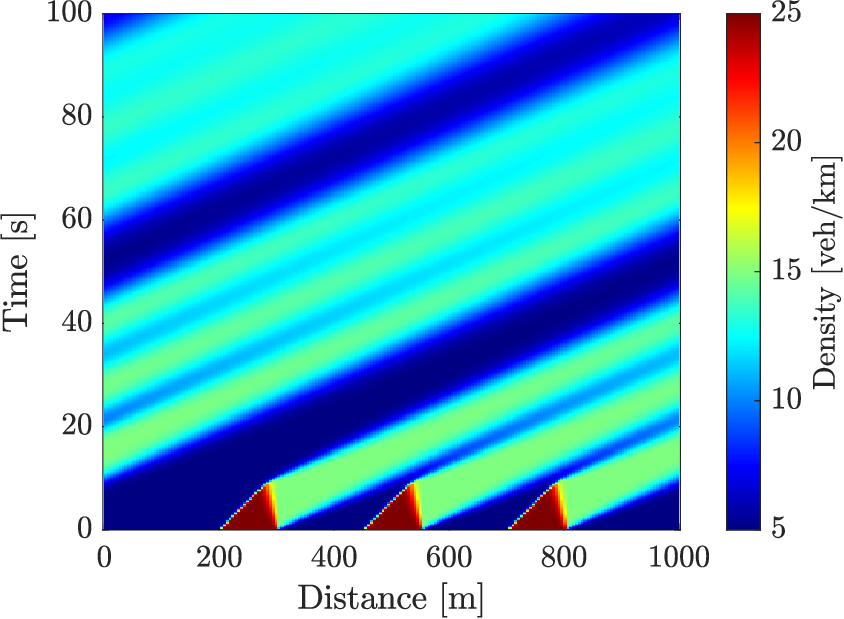}
            \caption{}
            \label{fig:CircularRdConRhoUnc}
    \end{subfigure}
    \begin{subfigure}{0.32\textwidth}
            \centering
            \includegraphics[width=\textwidth]{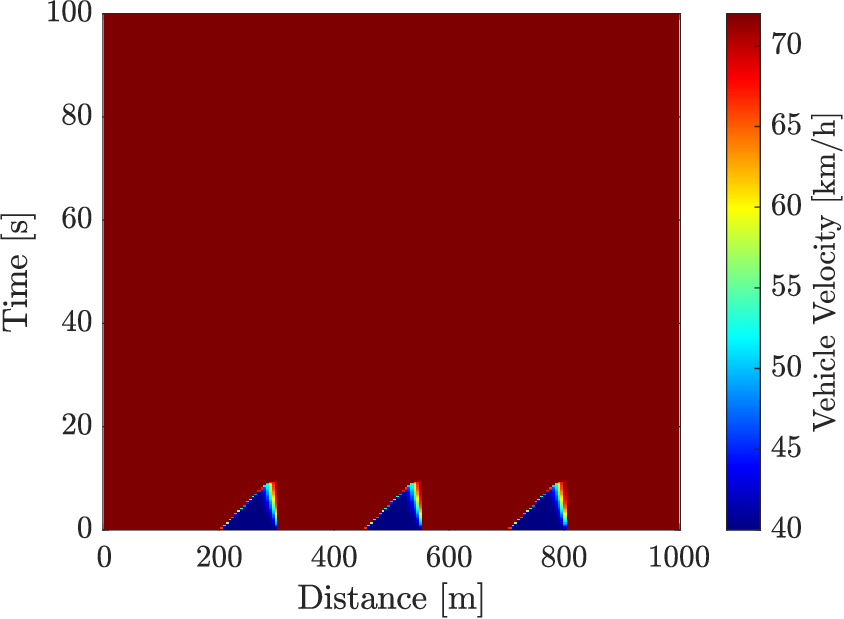}
            \caption{}
            \label{fig:CircularRdConVelUnc}
    \end{subfigure}
    \begin{subfigure}{0.32\textwidth}
            \centering
            \includegraphics[width=\textwidth]{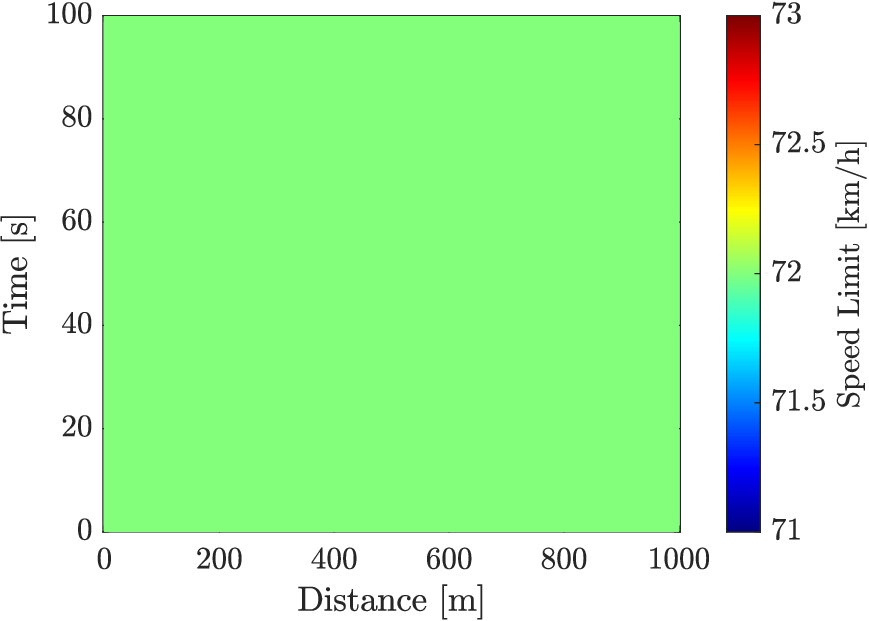}
            \caption{}
            \label{fig:CircularRdConLimUnc}
    \end{subfigure}
    \caption{Uncontrolled density (a), vehicle speed (b), and speed limit (c) for circular road case study.}
    \label{fig:CircularRdUnc}
\end{figure*}

\begin{figure*}[ht!]
    \centering
    \begin{subfigure}{0.32\textwidth}
            \centering
            \includegraphics[width=\textwidth]{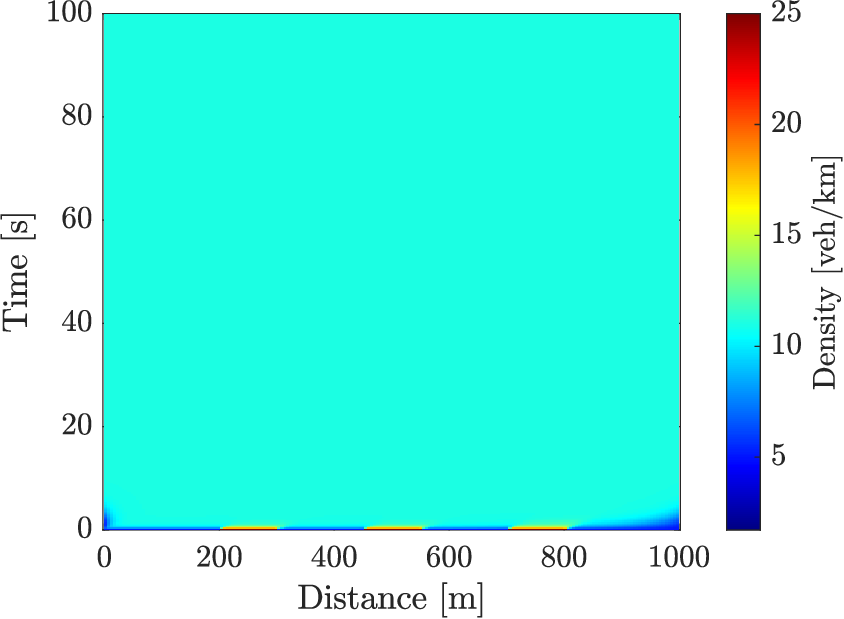}
            \caption{}
            \label{fig:CircularRdConRhoQ1e_2R1e_1}
    \end{subfigure}
    \begin{subfigure}{0.32\textwidth}
            \centering
            \includegraphics[width=\textwidth]{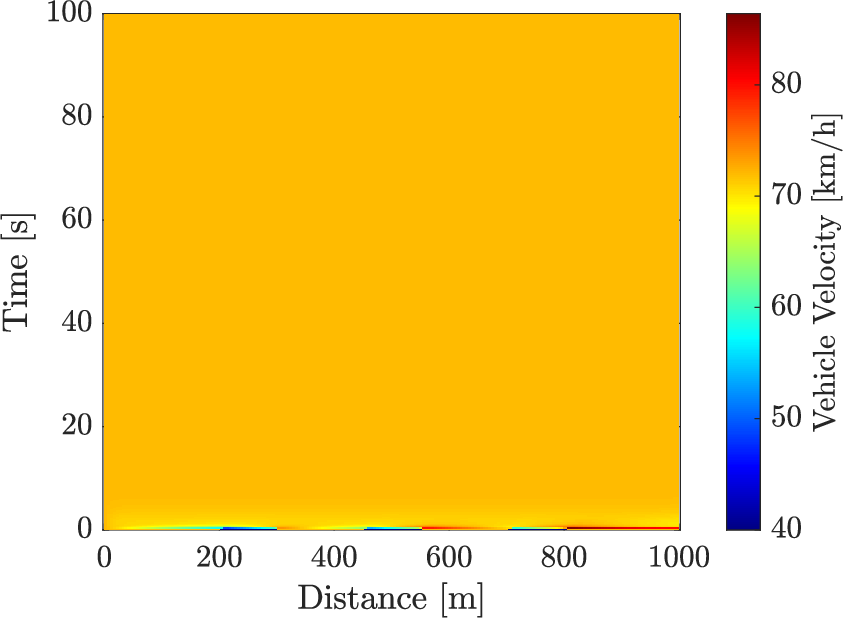}
            \caption{}
            \label{fig:CircularRdConVelQ1e_2R1e_1}
    \end{subfigure}
    \begin{subfigure}{0.32\textwidth}
            \centering
            \includegraphics[width=\textwidth]{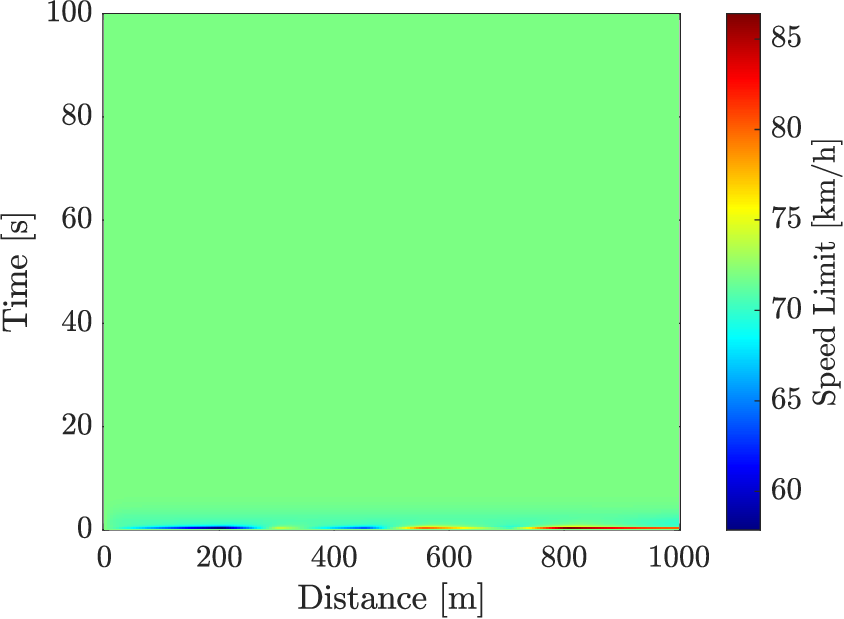}
            \caption{}
            \label{fig:CircularRdConLimQ1e_2R1e_1}
    \end{subfigure}
    \caption{Density (a), vehicle speed (b), and speed limit (c) for circular road case study with $Q = 1\cdot10^{-2}$ and $R = 1\cdot10^{-1}$.}
    \label{fig:CircularRdQ1e_2R1e_1}
\end{figure*}

\begin{figure*}[ht!]
    \centering
    \begin{subfigure}{0.32\textwidth}
            \centering
            \includegraphics[width=\textwidth]{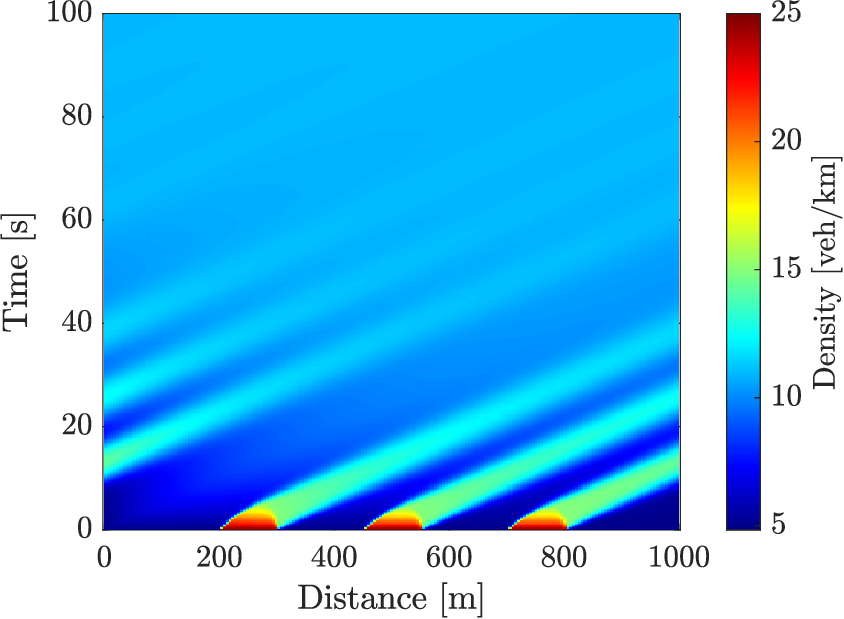}
            \caption{}
            \label{fig:CircularRdConRhoQ1e_4R1e_1}
    \end{subfigure}
    \begin{subfigure}{0.32\textwidth}
            \centering
            \includegraphics[width=\textwidth]{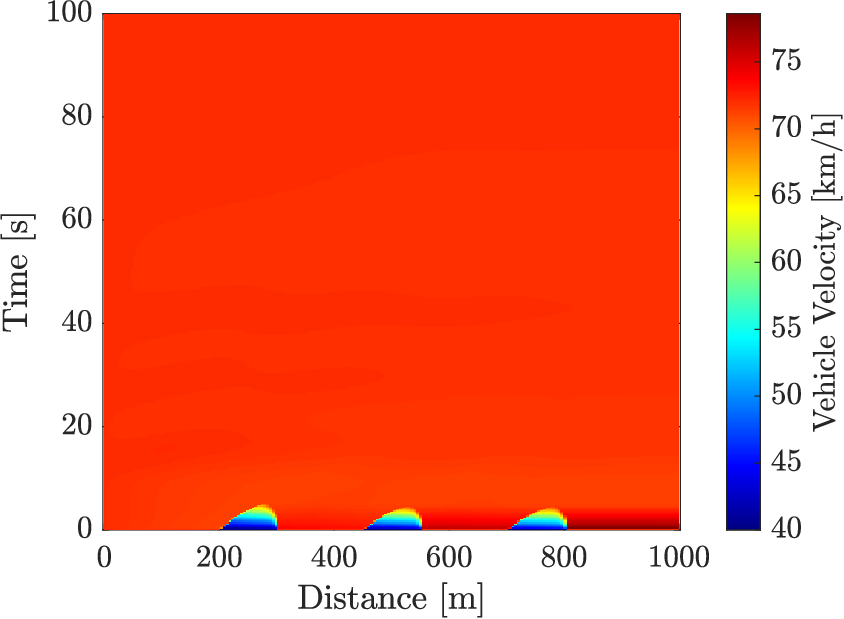}
            \caption{}
            \label{fig:CircularRdConVelQ1e_4R1e_1}
    \end{subfigure}
    \begin{subfigure}{0.32\textwidth}
            \centering
            \includegraphics[width=\textwidth]{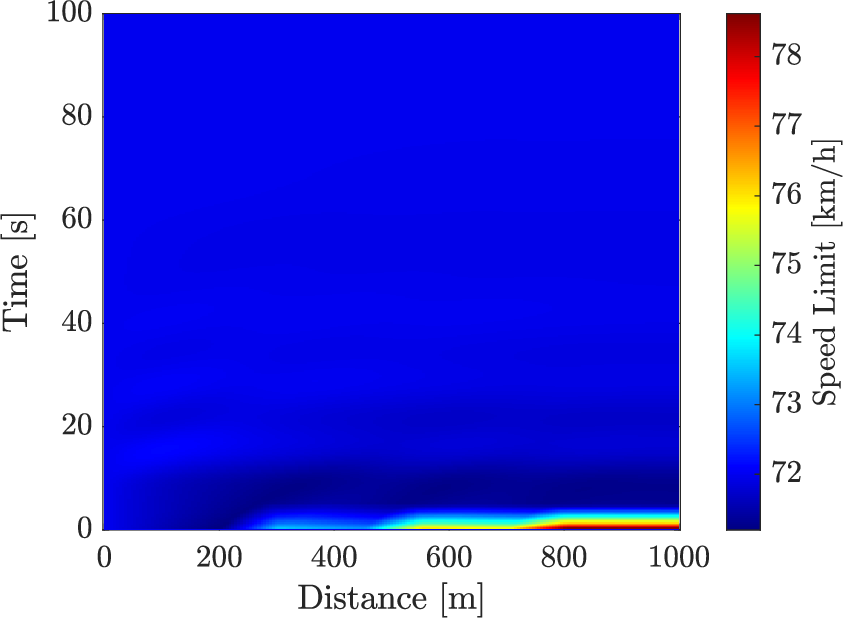}
            \caption{}
            \label{fig:CircularRdConLimQ1e_4R1e_1}
    \end{subfigure}
    \caption{Density (a), vehicle speed (b), and speed limit (c) for circular road case study with $Q = 1\cdot10^{-4}$ and $R = 1\cdot10^{-1}$.}
    \label{fig:CircularRdQ1e_4R1e_1}
\end{figure*}

\begin{figure*}[ht!]
    \centering
    \begin{subfigure}{0.32\textwidth}
            \centering
            \includegraphics[width=\textwidth]{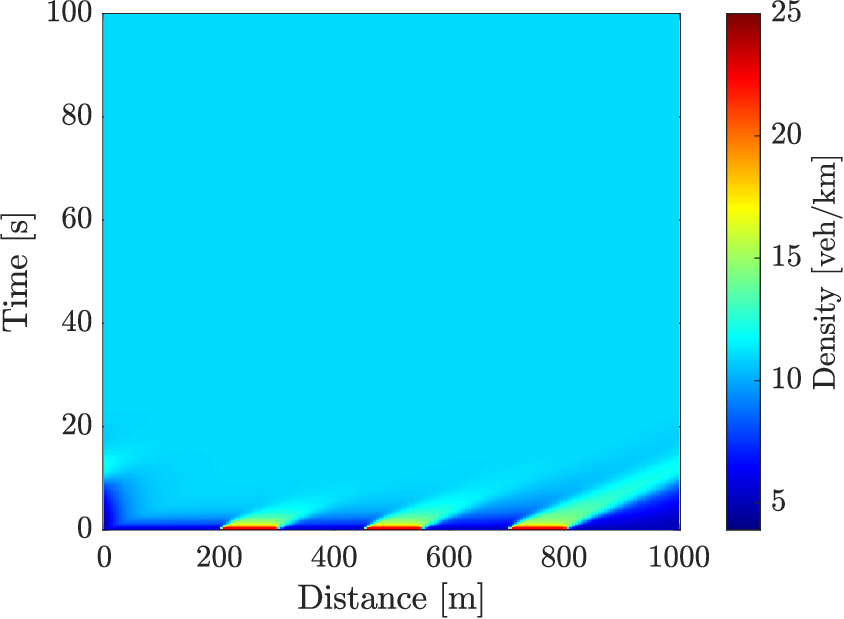}
            \caption{}
            \label{fig:CircularRdConRhoQ1e_4R1e_2}
    \end{subfigure}
    \begin{subfigure}{0.32\textwidth}
            \centering
            \includegraphics[width=\textwidth]{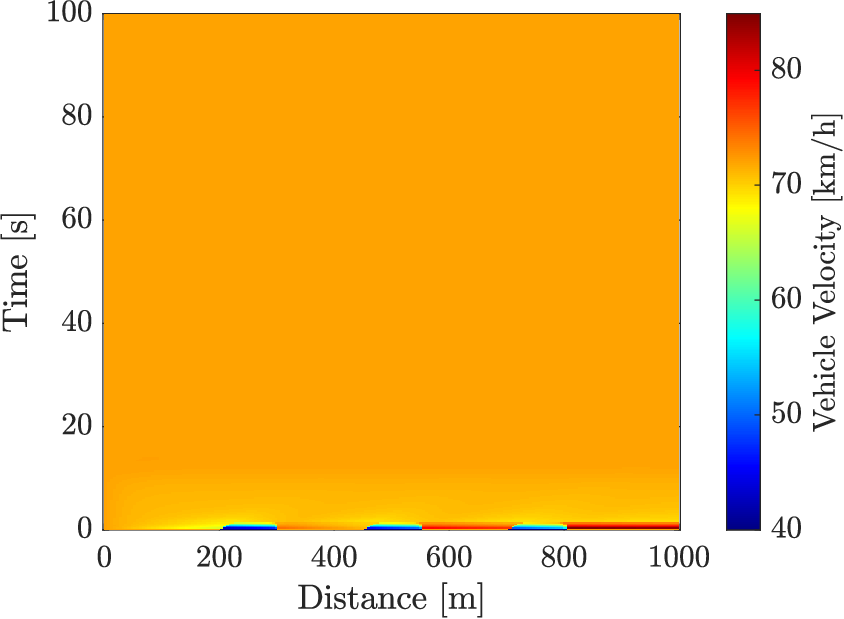}
            \caption{}
            \label{fig:CircularRdConVelQ1e_4R1e_2}
    \end{subfigure}
    \begin{subfigure}{0.32\textwidth}
            \centering
            \includegraphics[width=\textwidth]{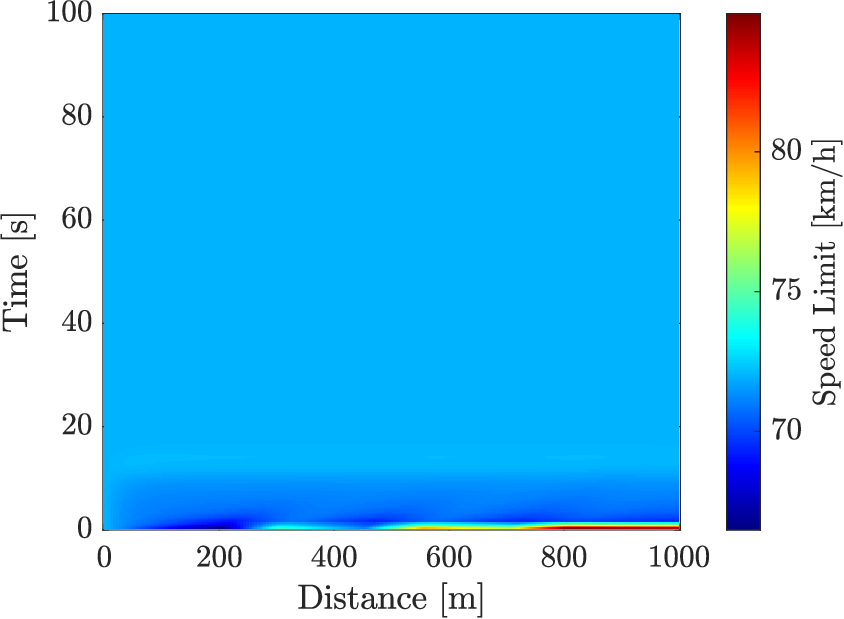}
            \caption{}
            \label{fig:CircularRdConLimQ1e_4R1e_2}
    \end{subfigure}
    \caption{Density (a), vehicle speed (b), and speed limit (c) for circular road case study with $Q = 1\cdot10^{-4}$ and $R = 1\cdot10^{-2}$.}
    \label{fig:CircularRdQ1e_4R1e_2}
\end{figure*}

The uncontrolled density and vehicle speed for the circular road case study is shown in \cref{fig:CircularRdUnc}. As can be seen, the high density regions cause waves of traffic at a higher density than the desired density. \cref{fig:CircularRdQ1e_2R1e_1}-\cref{fig:CircularRdQ1e_4R1e_1} show three different cases with varying control aggressiveness. In \cref{fig:CircularRdQ1e_2R1e_1} the controller is too aggressive as the penalty on state deviation is high and the use of the control action is not penalized enough. The maximum speed limit, shown in \cref{fig:CircularRdConLimQ1e_2R1e_1}, also grows too high as it exceeds 85 km/h when the nominal is 72 km/h. \cref{fig:CircularRdQ1e_4R1e_1} presents a much more reasonable case, where the maximum speed limit does not exceed 80 km/h and it is still able to regulate the density to the desired setpoint. The penalty on the state deviation is kept the same in \cref{fig:CircularRdQ1e_4R1e_2}, but the penalty on the control input is relaxed. This results in, again, a higher maximum speed limit as seen in \cref{fig:CircularRdConLimQ1e_4R1e_2}.

\begin{figure}[ht!]
    \centering
    \begin{subfigure}{0.45\columnwidth}
            \centering
            \includegraphics[width=\textwidth]{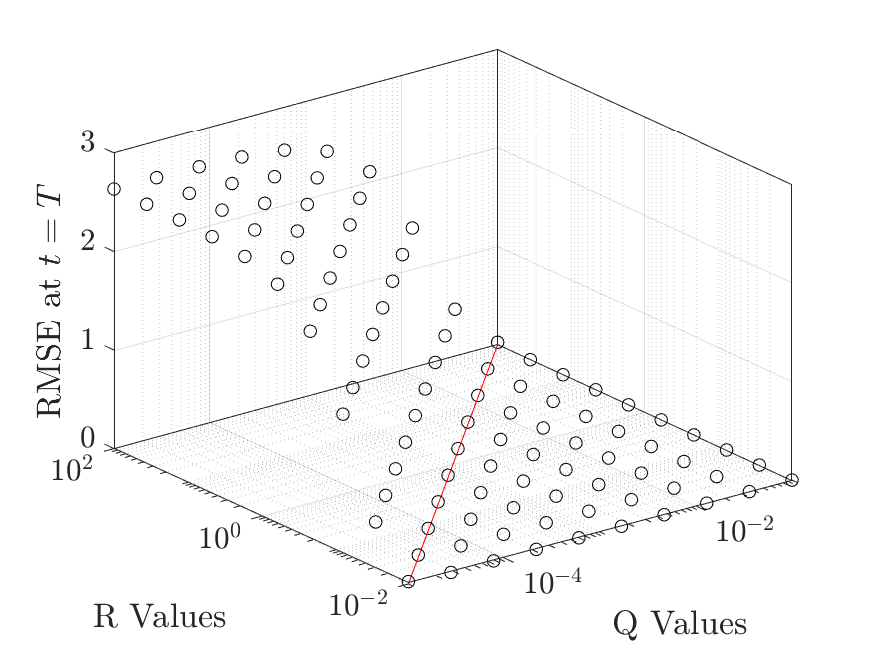}
            \caption{}
            \label{fig:FinalTimeRMSECircular}
    \end{subfigure}
    \begin{subfigure}{0.45\columnwidth}
            \centering
            \includegraphics[width=\textwidth]{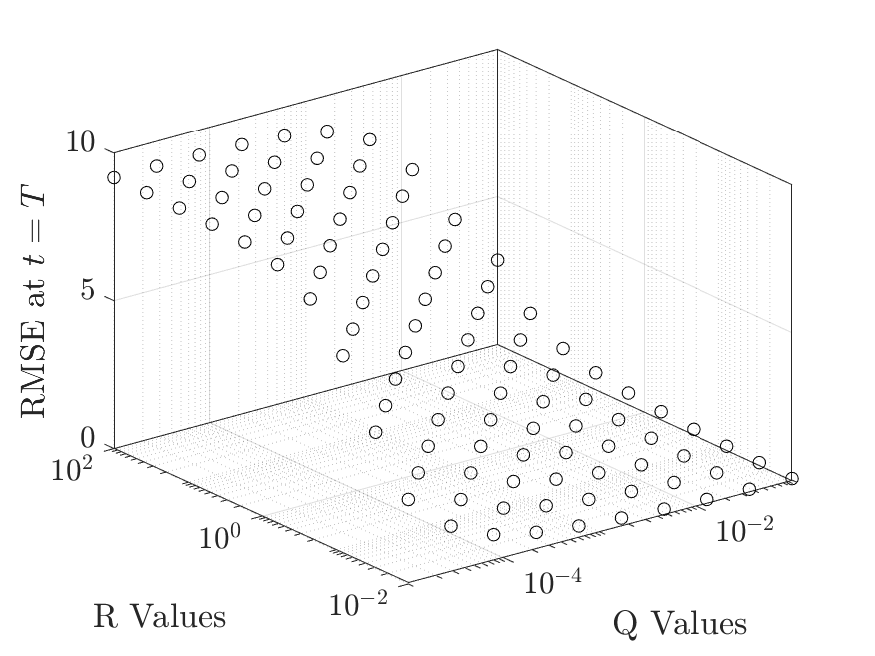}
            \caption{}
            \label{fig:FinalTimeRMSEStraight}
    \end{subfigure}
    \caption{RMSE at the final time step for a range of $Q$ and $R$ values for circular road (a) and straight road (b). }
    \label{fig:RMSEplots}
\end{figure}

The efficacy of different values of $Q$ and $R$ on the control goal is shown in \cref{fig:FinalTimeRMSECircular} by the the RMSE between the final density profile and the desired density profile . In \cref{fig:FinalTimeRMSECircular} the red line depicts the line along which the smallest control action, i.e. lowest maximum speed limit, is used while still achieving the goal in the simulation time.

\subsubsection{Straight Road Case Study}

\begin{figure*}[ht!]
    \centering
    \begin{subfigure}{0.32\textwidth}
            \centering
            \includegraphics[width=\textwidth]{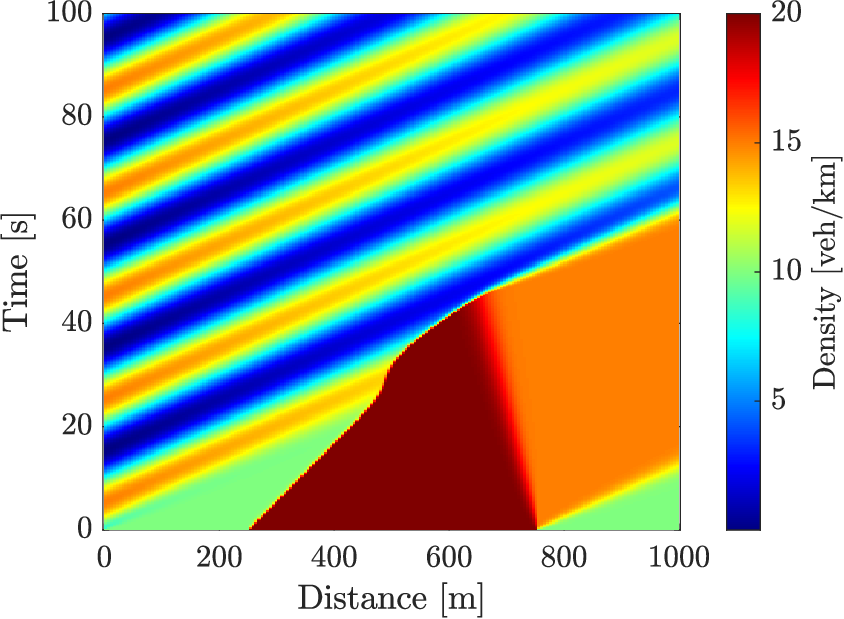}
            \caption{}
            \label{fig:StraightRdConRhoUnc}
    \end{subfigure}
    \begin{subfigure}{0.32\textwidth}
            \centering
            \includegraphics[width=\textwidth]{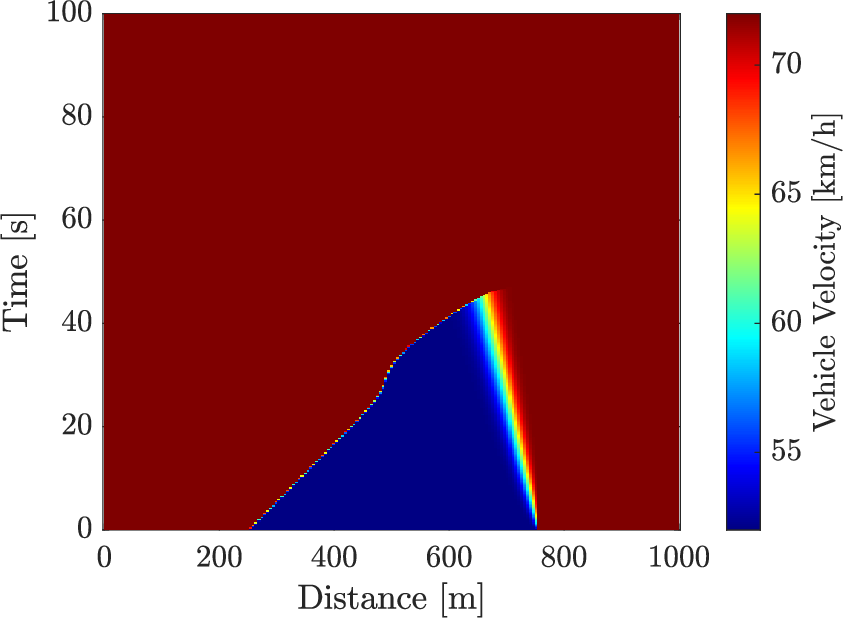}
            \caption{}
            \label{fig:StraightRdConVelUnc}
    \end{subfigure}
    \begin{subfigure}{0.32\textwidth}
            \centering
            \includegraphics[width=\textwidth]{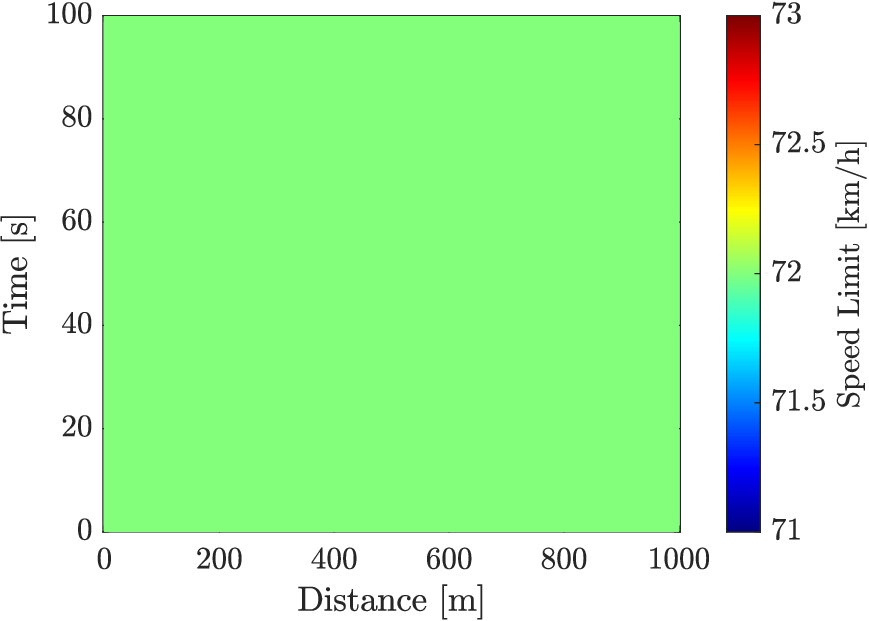}
            \caption{}
            \label{fig:StraightRdConLimUnc}
    \end{subfigure}
    \caption{Uncontrolled density (a), vehicle speed (b), and speed limit (c) for straight road case study.}
    \label{fig:StraightRdUnc}
\end{figure*}

\begin{figure*}[ht!]
    \centering
    \begin{subfigure}{0.32\textwidth}
            \centering
            \includegraphics[width=\textwidth]{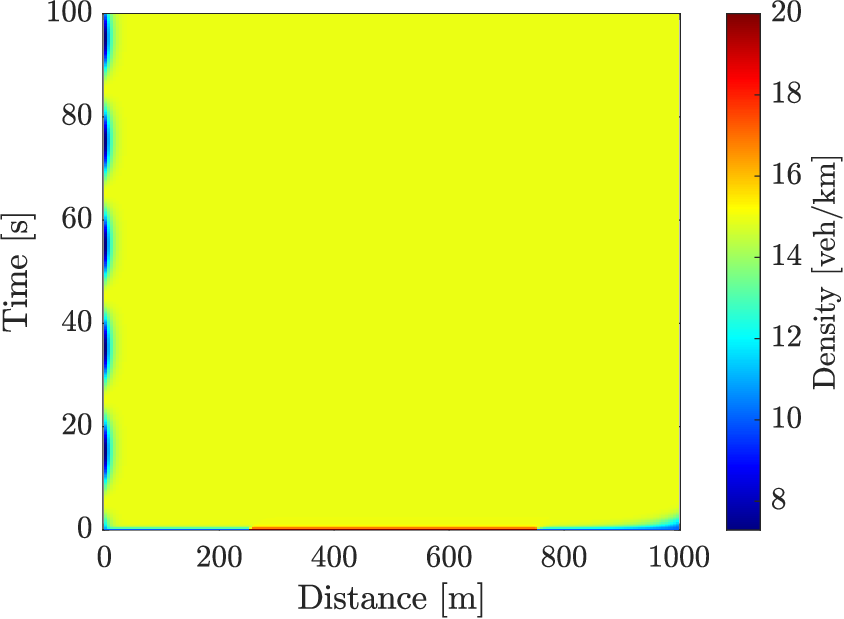}
            \caption{}
            \label{fig:StraightRdConRhoQ1e_2R1e_1}
    \end{subfigure}
    \begin{subfigure}{0.32\textwidth}
            \centering
            \includegraphics[width=\textwidth]{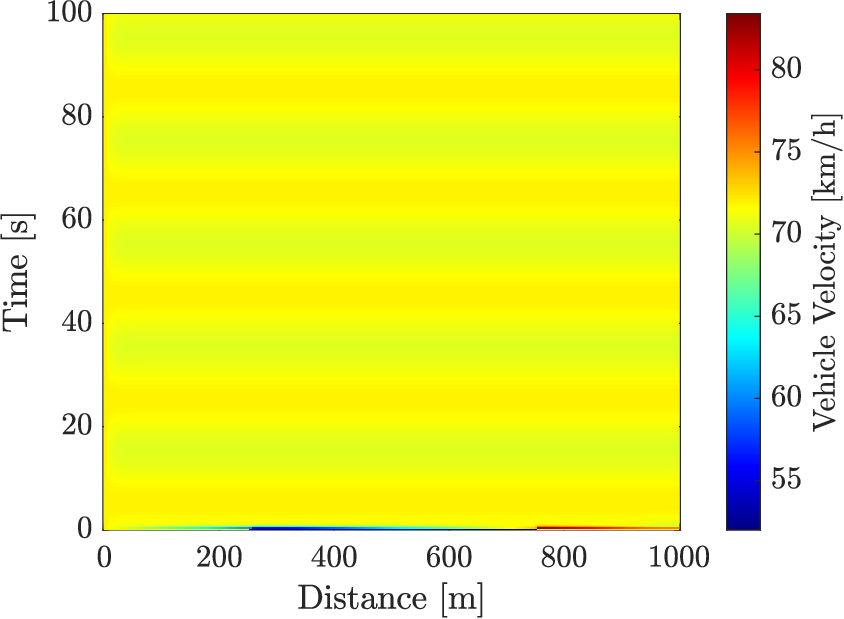}
            \caption{}
            \label{fig:StraightRdConVelQ1e_2R1e_1}
    \end{subfigure}
    \begin{subfigure}{0.32\textwidth}
            \centering
            \includegraphics[width=\textwidth]{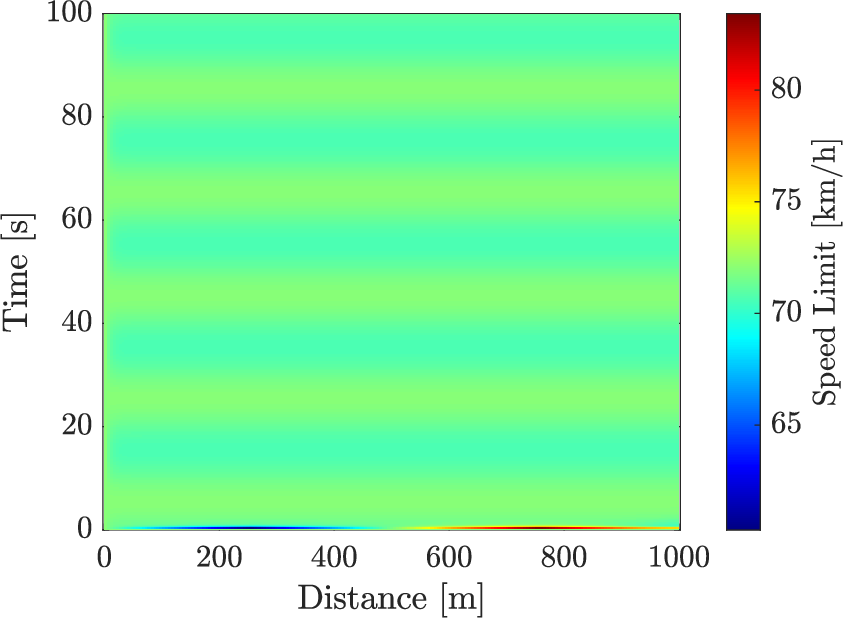}
            \caption{}
            \label{fig:StraightRdConLimQ1e_2R1e_1}
    \end{subfigure}
    \caption{Density (a), vehicle speed (b), and speed limit (c) for straight road case study with $Q = 1\cdot10^{-2}$ and $R = 1\cdot10^{-1}$.}
    \label{fig:StraightRdQ1e_2R1e_1}
\end{figure*}

\begin{figure*}[ht!]
    \centering
    \begin{subfigure}{0.32\textwidth}
            \centering
            \includegraphics[width=\textwidth]{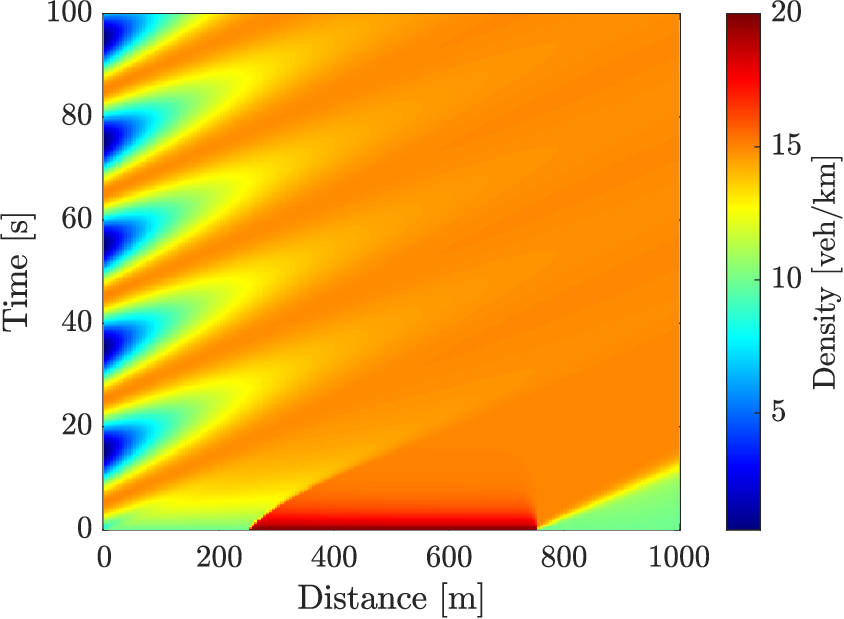}
            \caption{}
            \label{fig:StraightRdConRhoQ1e_4R1e_1}
    \end{subfigure}
    \begin{subfigure}{0.32\textwidth}
            \centering
            \includegraphics[width=\textwidth]{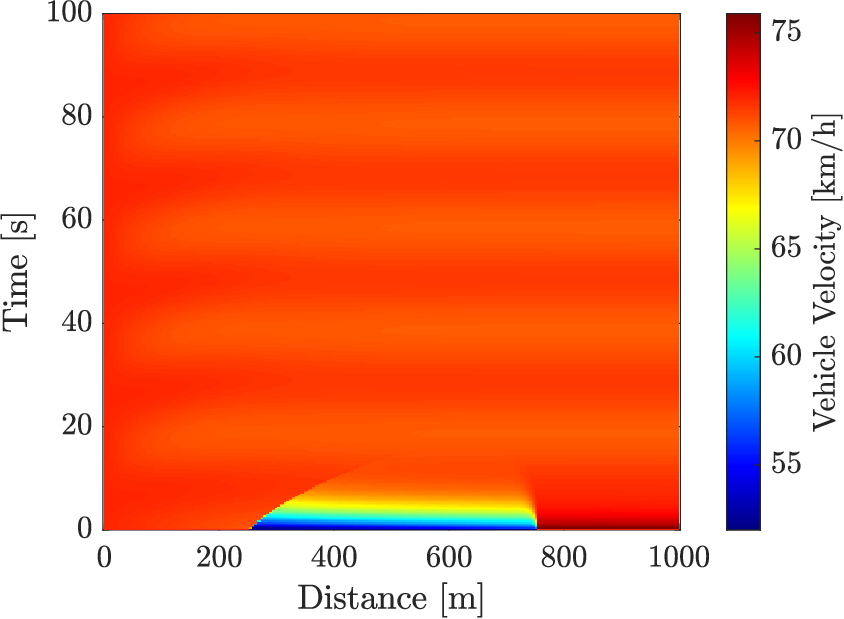}
            \caption{}
            \label{fig:StraightRdConVelQ1e_4R1e_1}
    \end{subfigure}
    \begin{subfigure}{0.32\textwidth}
            \centering
            \includegraphics[width=\textwidth]{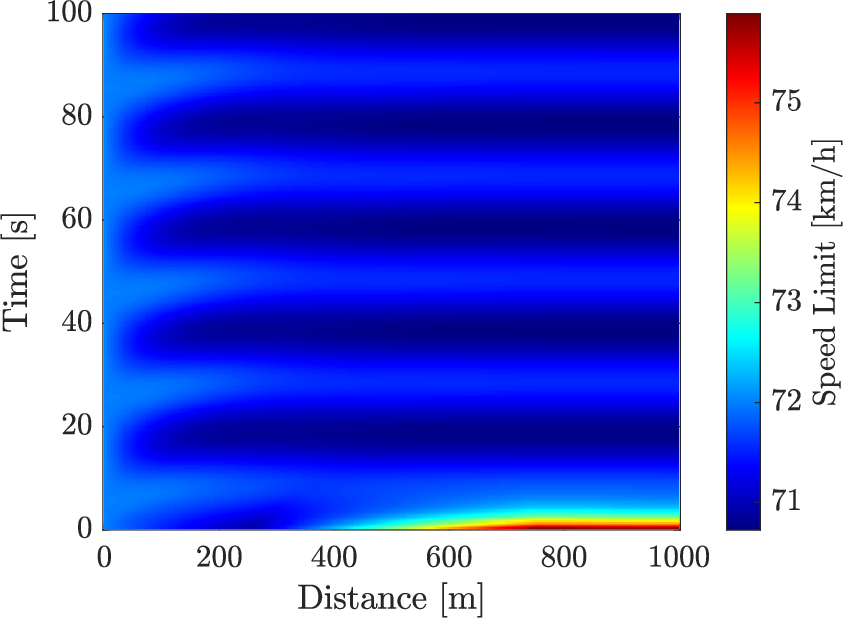}
            \caption{}
            \label{fig:StraightRdConLimQ1e_4R1e_1}
    \end{subfigure}
    \caption{Density (a), vehicle speed (b), and speed limit (c) for straight road case study with $Q = 1\cdot10^{-4}$ and $R = 1\cdot10^{-1}$.}
    \label{fig:StraightRdQ1e_4R1e_1}
\end{figure*}

\begin{figure*}[ht!]
    \centering
    \begin{subfigure}{0.32\textwidth}
            \centering
            \includegraphics[width=\textwidth]{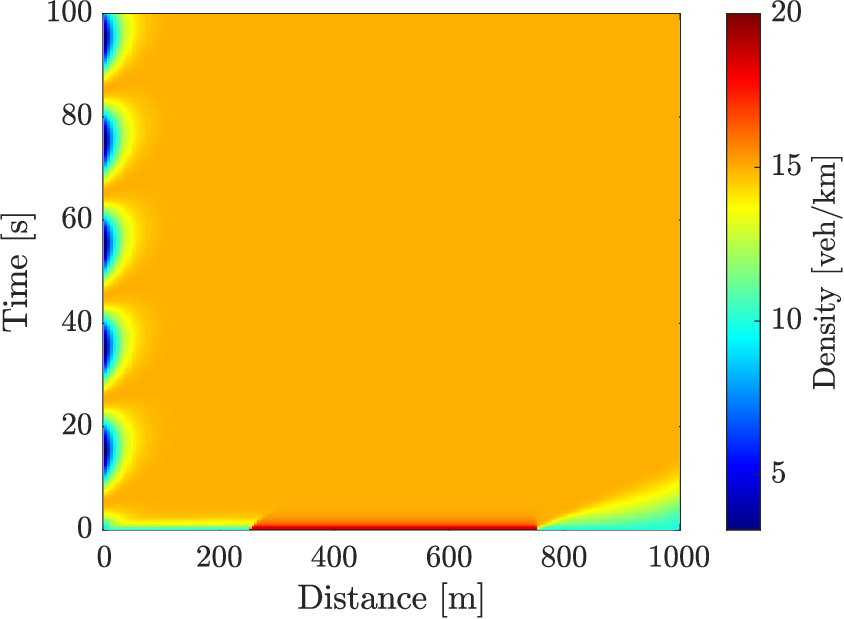}
            \caption{}
            \label{fig:StraightRdConRhoQ1e_4R1e_2}
    \end{subfigure}
    \begin{subfigure}{0.32\textwidth}
            \centering
            \includegraphics[width=\textwidth]{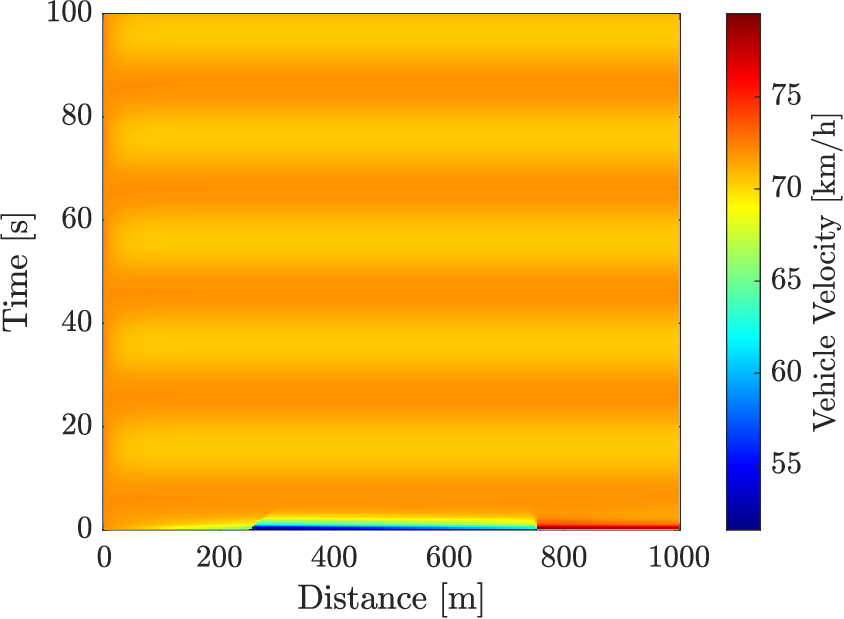}
            \caption{}
            \label{fig:StraightRdConVelQ1e_4R1e_2}
    \end{subfigure}
    \begin{subfigure}{0.32\textwidth}
            \centering
            \includegraphics[width=\textwidth]{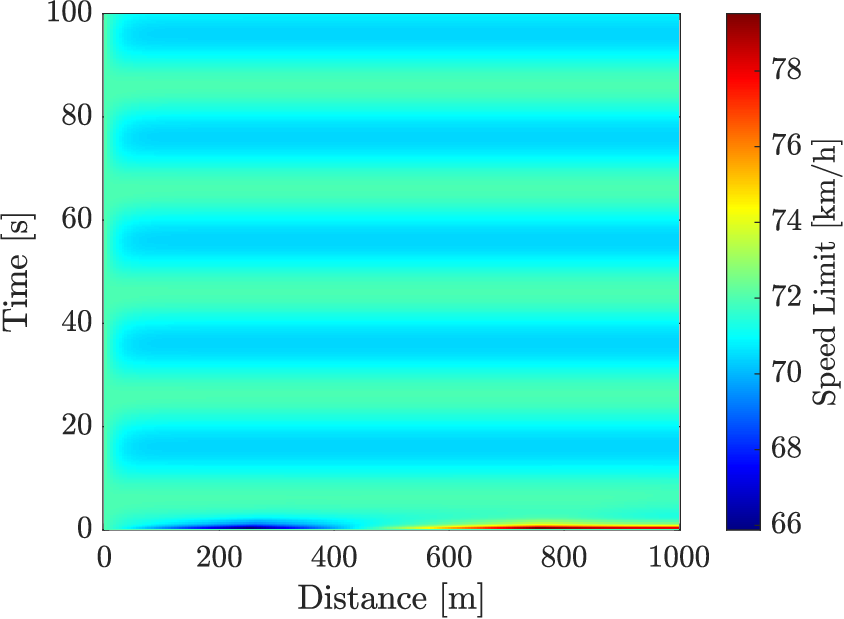}
            \caption{}
            \label{fig:StraightRdConLimQ1e_4R1e_2}
    \end{subfigure}
    \caption{Density (a), vehicle speed (b), and speed limit (c) for straight road case study with $Q = 1\cdot10^{-4}$ and $R = 1\cdot10^{-2}$.}
    \label{fig:StraightRdQ1e_4R1e_2}
\end{figure*}

In the straight road case, not only does the controller have to deal with a high density region in the domain, but also with incoming traffic at the boundary. The uncontrolled scenario is shown in \cref{fig:StraightRdUnc} where the incoming flow catches up and meets the high-density traffic already in the domain. In \cref{fig:StraightRdQ1e_2R1e_1}, a very aggressive controller is shown that alleviates the high density traffic almost immediately and greatly reduces the impact of the incoming boundary flow. A more reasonable controller is shown in \cref{fig:StraightRdQ1e_4R1e_1}, which only increases the speed limit by around 4 km/h. While this same control action worked well for the circular road, it is less effective in the straight-road case because it cannot dissipate the continuous traffic inflow fast enough. One interesting aspect to note in \cref{fig:StraightRdConLimQ1e_4R1e_1} is that there are bands of different speed limits that change every 15 to 20 s. As with the circular road, we then relax the penalty on the control input and, as seen in \cref{fig:StraightRdQ1e_4R1e_2}, see that the control action grows. Here, though, it still commands reasonable speed limits while still mitigating the impact of the incoming flow.

An interesting difference with this case study from the circular road case study is seen in \cref{fig:FinalTimeRMSEStraight}. Before, there was a clear line delineating the least aggressive control action that still achieves the control goal, here all but the strongest control actions fail to completely regulate the traffic to the desired density since he controller does not have knowledge of the incoming traffic, and thus cannot account for it.


\subsection{Impact of $t_f$ and $S$ on Finite Horizon Optimal Solution}

Since the impact of $Q$ and $R$ on the feedback solution was investigated for the infinite horizon case, here the focus is on the final time and the value of $S$. Since $Q=1\cdot10^{-4}$ and $R=1\cdot10^{-1}$ gave reasonable results for the infinite horizon controller, they are used here. The initial and boundary conditions are kept the same as before, as well as the model parameters. To see the impact of final time on the finite horizon controller, four different final times $t_f\in[100 \; 75 \; 50 \; 25]$ s are used keeping $S=R$ in all cases. Then, to investigate the impact of $S$ on the solution, $S$ is changed to $S=2R$ and $S=0.5R$ for the case when $t_f = 50$ s.

\subsubsection{Circular Road Case Study}

\begin{figure*}[ht!]
    \centering
    \begin{subfigure}{0.32\textwidth}
            \centering
            \includegraphics[width=\textwidth]{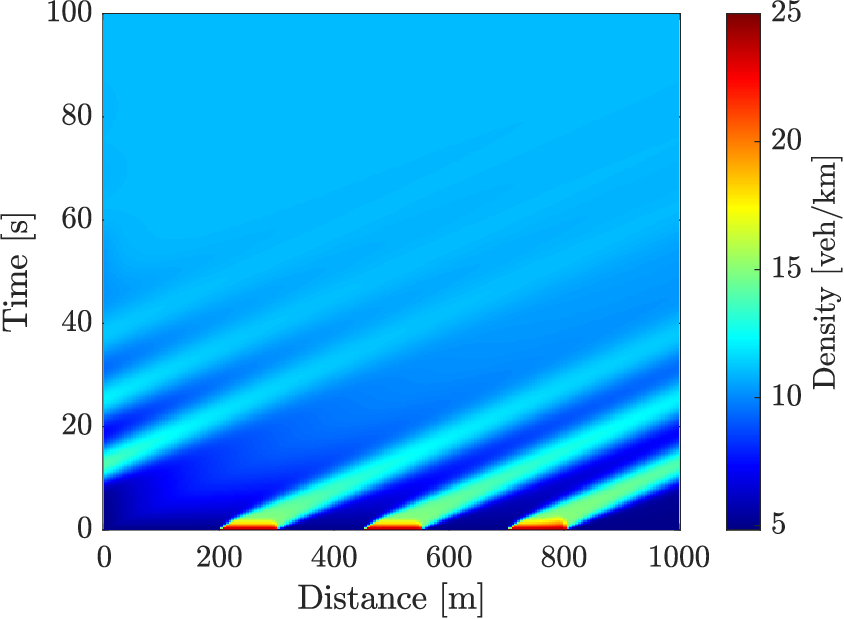}
            \caption{}
            \label{fig:FinCircularRdConRhot_f100SR}
    \end{subfigure}
    \begin{subfigure}{0.32\textwidth}
            \centering
            \includegraphics[width=\textwidth]{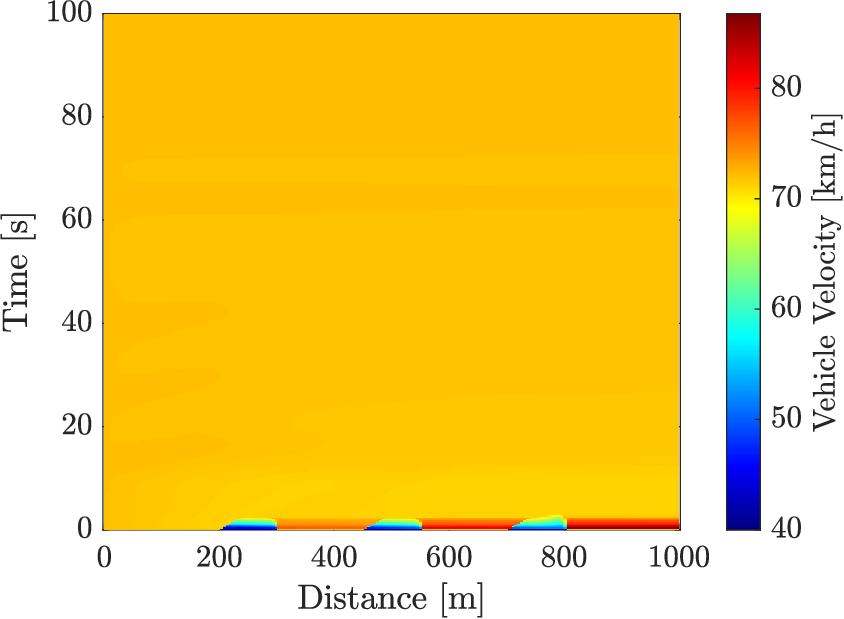}
            \caption{}
            \label{fig:FinCircularRdConVelt_f100SR}
    \end{subfigure}
    \begin{subfigure}{0.32\textwidth}
            \centering
            \includegraphics[width=\textwidth]{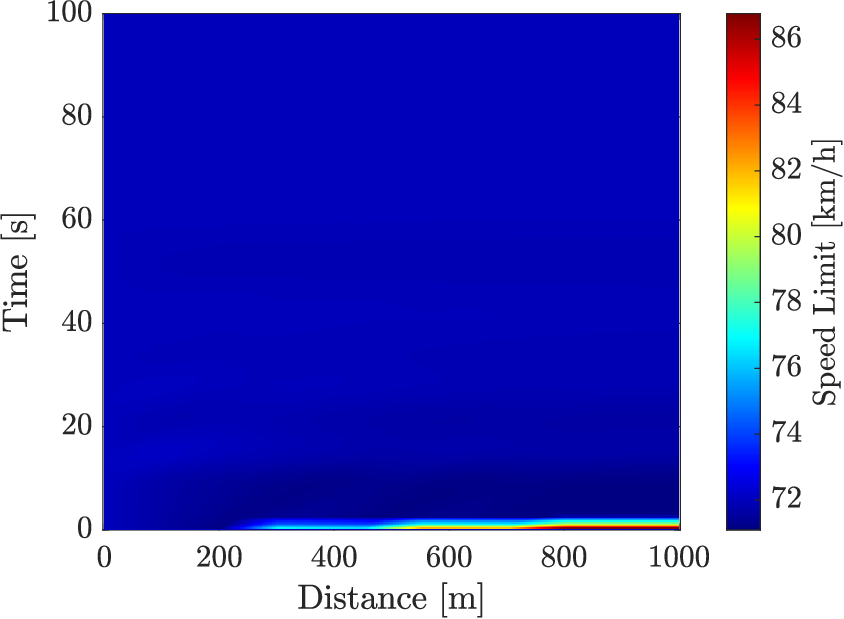}
            \caption{}
            \label{fig:FinCircularRdConLimt_f100SR}
    \end{subfigure}
    \caption{Density (a), vehicle speed (b), and speed limit (c) for finite-horizon control on a circular road with $t_f = 100$ s, $S = R$.}
    \label{fig:FinCircularRdt_f100SR}
\end{figure*}

\begin{figure*}[ht!]
    \centering
    \begin{subfigure}{0.32\textwidth}
            \centering
            \includegraphics[width=\textwidth]{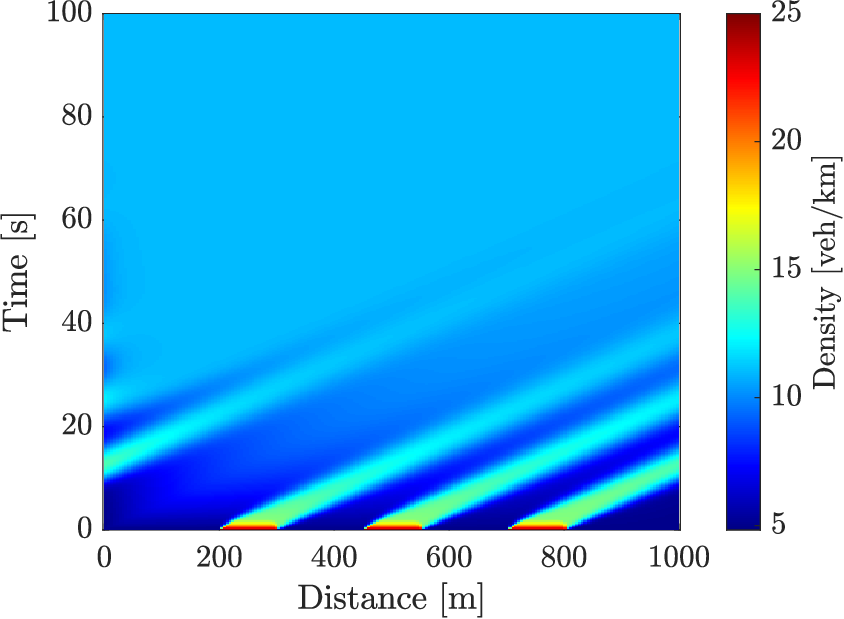}
            \caption{}
            \label{fig:FinCircularRdConRhot_f75SR}
    \end{subfigure}
    \begin{subfigure}{0.32\textwidth}
            \centering
            \includegraphics[width=\textwidth]{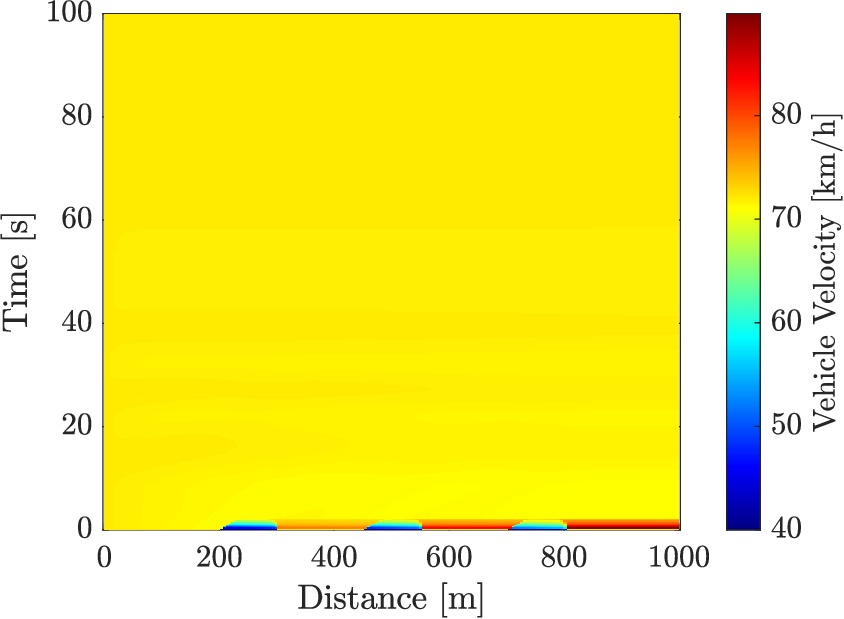}
            \caption{}
            \label{fig:FinCircularRdConVelt_f75SR}
    \end{subfigure}
    \begin{subfigure}{0.32\textwidth}
            \centering
            \includegraphics[width=\textwidth]{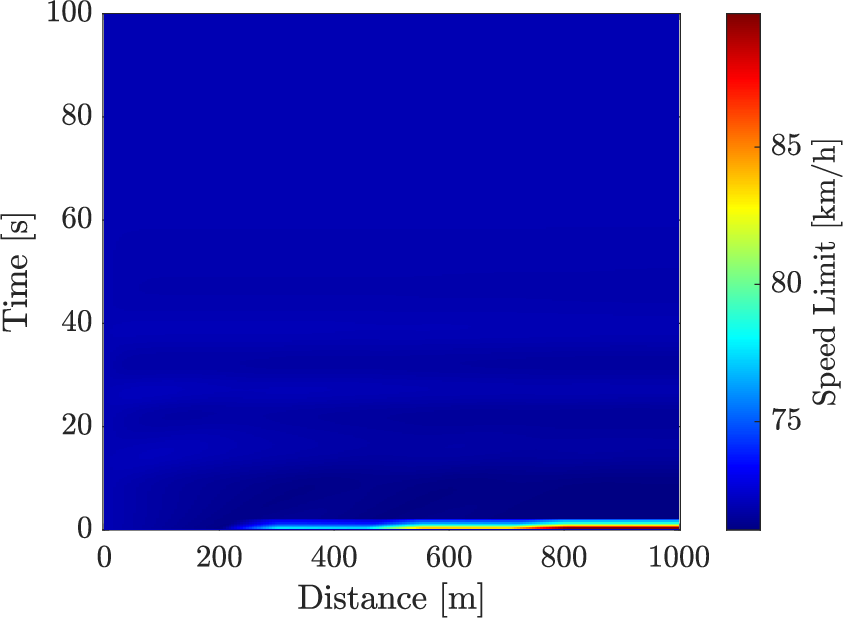}
            \caption{}
            \label{fig:FinCircularRdConLimt_f75SR}
    \end{subfigure}
    \caption{Density (a), vehicle speed (b), and speed limit (c) for finite-horizon control on a circular road with $t_f = 75$ s, $S = R$.}
    \label{fig:FinCircularRdt_f75SR}
\end{figure*}

\begin{figure*}[ht!]
    \centering
    \begin{subfigure}{0.32\textwidth}
            \centering
            \includegraphics[width=\textwidth]{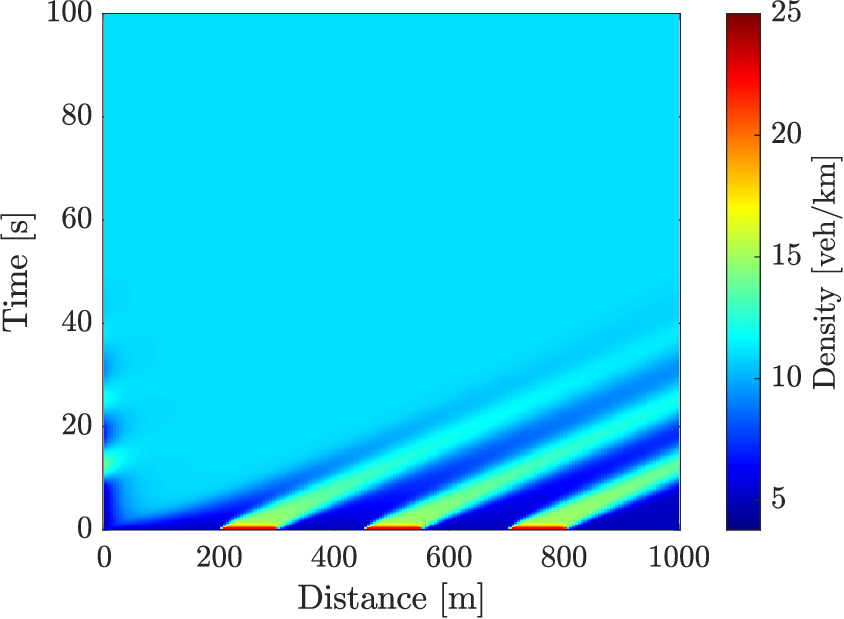}
            \caption{}
            \label{fig:FinCircularRdConRhot_f50SR}
    \end{subfigure}
    \begin{subfigure}{0.32\textwidth}
            \centering
            \includegraphics[width=\textwidth]{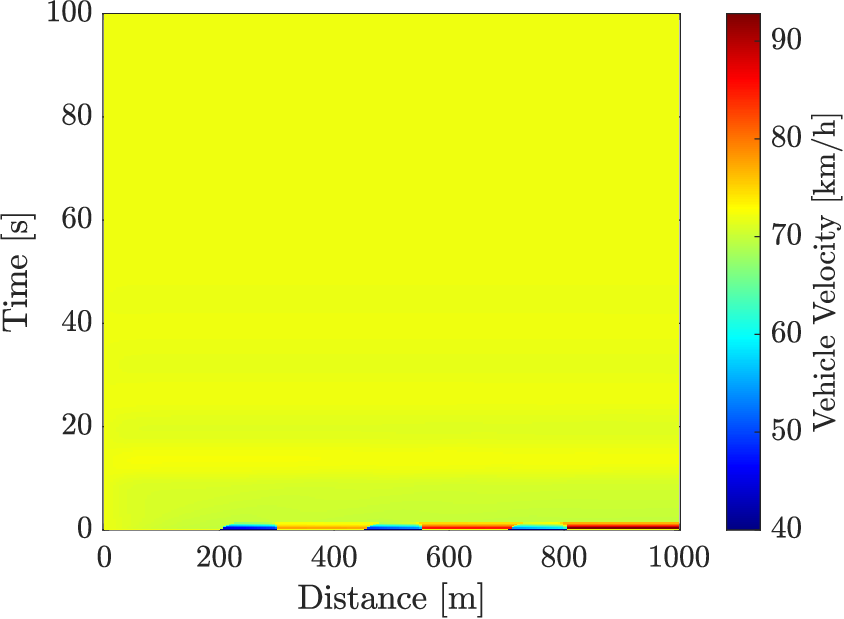}
            \caption{}
            \label{fig:FinCircularRdConVelt_f50SR}
    \end{subfigure}
    \begin{subfigure}{0.32\textwidth}
            \centering
            \includegraphics[width=\textwidth]{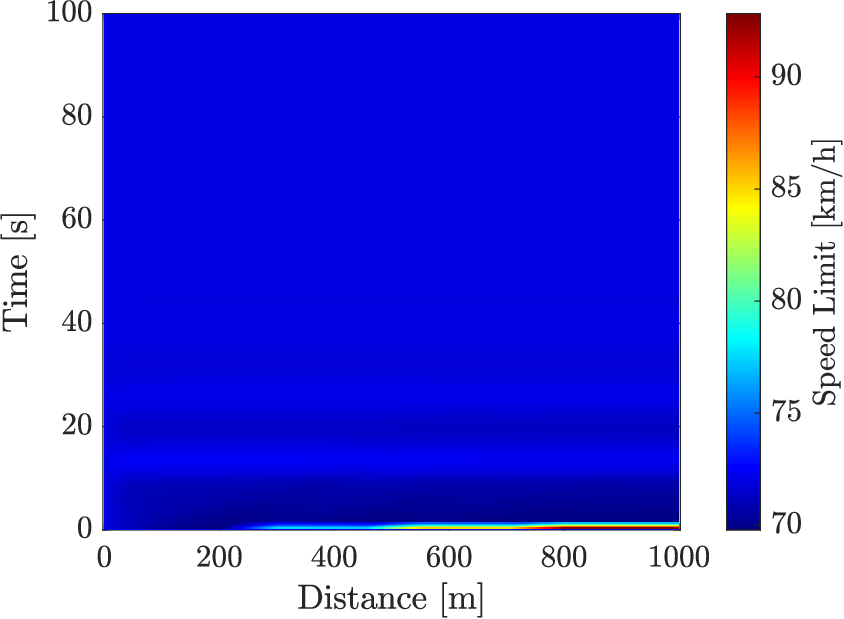}
            \caption{}
            \label{fig:FinCircularRdConLimt_f50SR}
    \end{subfigure}
    \caption{Density (a), vehicle speed (b), and speed limit (c) for finite-horizon control on a circular road with $t_f = 50$ s, $S = R$.}
    \label{fig:FinCircularRdt_f50SR}
\end{figure*}

\begin{figure*}[ht!]
    \centering
    \begin{subfigure}{0.32\textwidth}
            \centering
            \includegraphics[width=\textwidth]{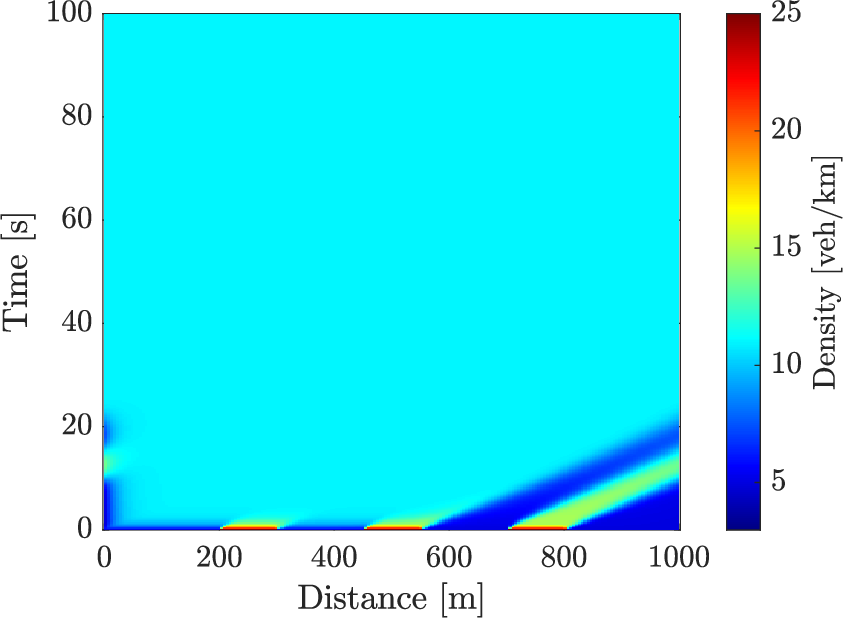}
            \caption{}
            \label{fig:FinCircularRdConRhot_f25SR}
    \end{subfigure}
    \begin{subfigure}{0.32\textwidth}
            \centering
            \includegraphics[width=\textwidth]{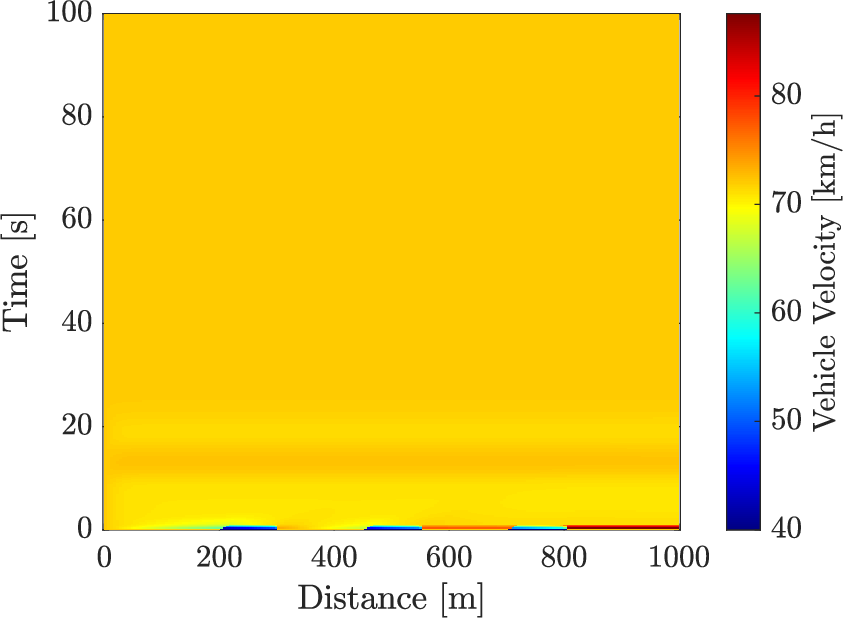}
            \caption{}
            \label{fig:FinCircularRdConVelt_f25SR}
    \end{subfigure}
    \begin{subfigure}{0.32\textwidth}
            \centering
            \includegraphics[width=\textwidth]{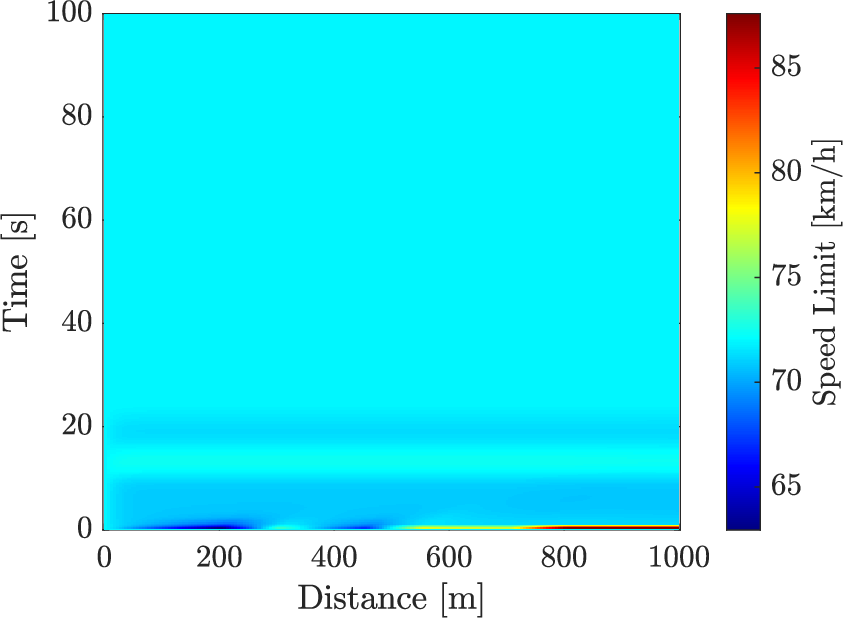}
            \caption{}
            \label{fig:FinCircularRdConLimt_f25SR}
    \end{subfigure}
    \caption{Density (a), vehicle speed (b), and speed limit (c) for finite-horizon control on a circular road with $t_f = 25$ s, $S = R$.}
    \label{fig:FinCircularRdt_f25SR}
\end{figure*}

\begin{figure*}[ht!]
    \centering
    \begin{subfigure}{0.32\textwidth}
            \centering
            \includegraphics[width=\textwidth]{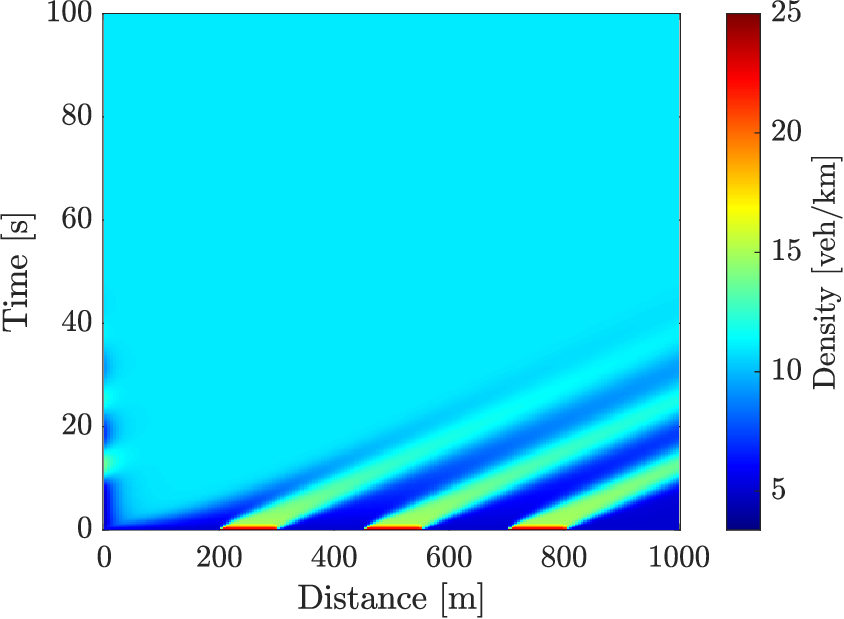}
            \caption{}
            \label{fig:FinCircularRdConRhot_f50S2R}
    \end{subfigure}
    \begin{subfigure}{0.32\textwidth}
            \centering
            \includegraphics[width=\textwidth]{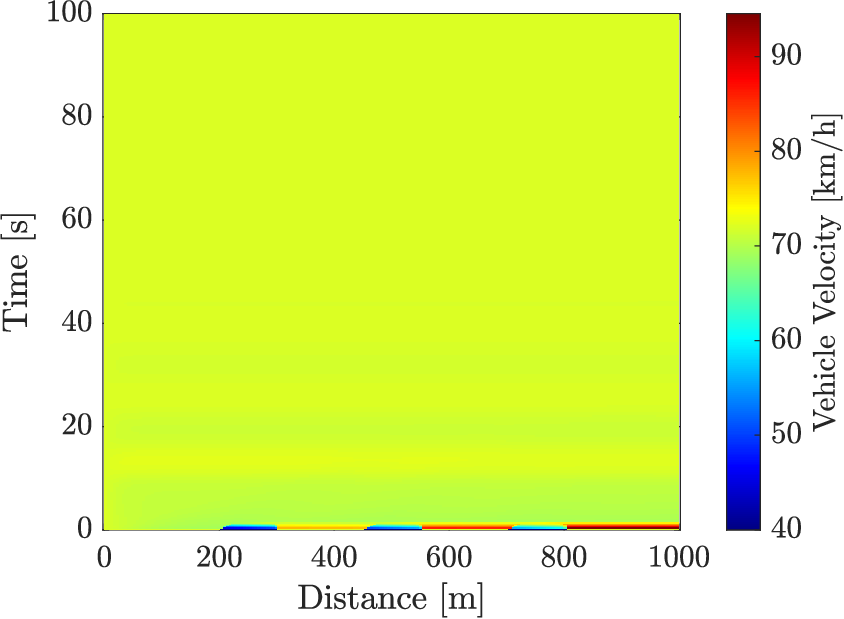}
            \caption{}
            \label{fig:FinCircularRdConVelt_f50S2R}
    \end{subfigure}
    \begin{subfigure}{0.32\textwidth}
            \centering
            \includegraphics[width=\textwidth]{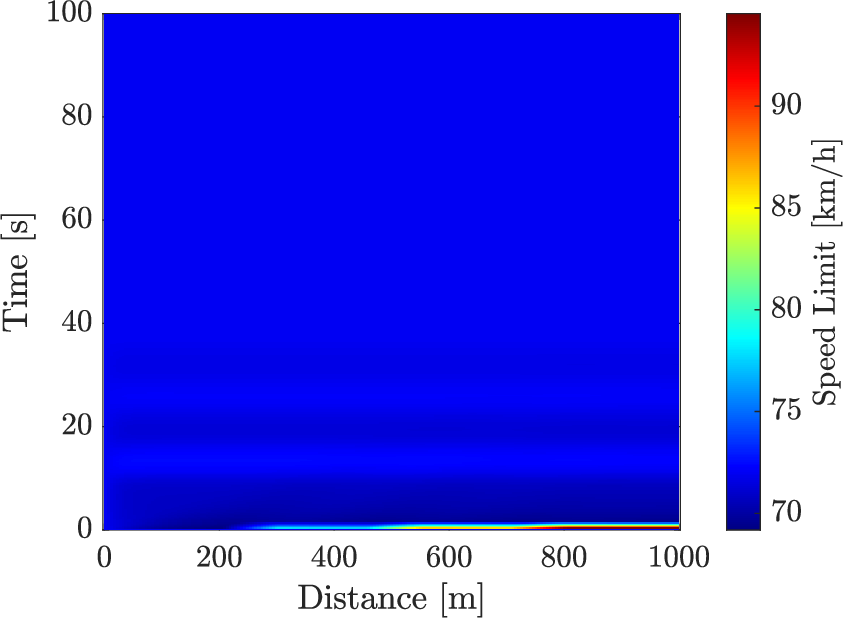}
            \caption{}
            \label{fig:FinCircularRdConLimt_f50S2R}
    \end{subfigure}
    \caption{Density (a), vehicle speed (b), and speed limit (c) for finite-horizon control on a circular road with $t_f = 50$ s, $S = 2R$.}
    \label{fig:FinCircularRdt_f50S2R}
\end{figure*}

\begin{figure*}[ht!]
    \centering
    \begin{subfigure}{0.32\textwidth}
            \centering
            \includegraphics[width=\textwidth]{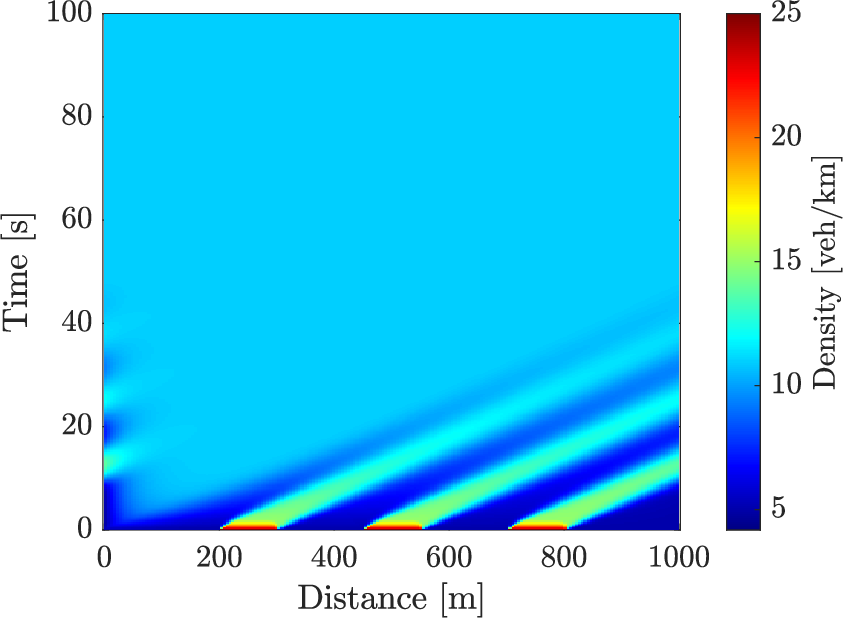}
            \caption{}
            \label{fig:FinCircularRdConRhot_f50S05R}
    \end{subfigure}
    \begin{subfigure}{0.32\textwidth}
            \centering
            \includegraphics[width=\textwidth]{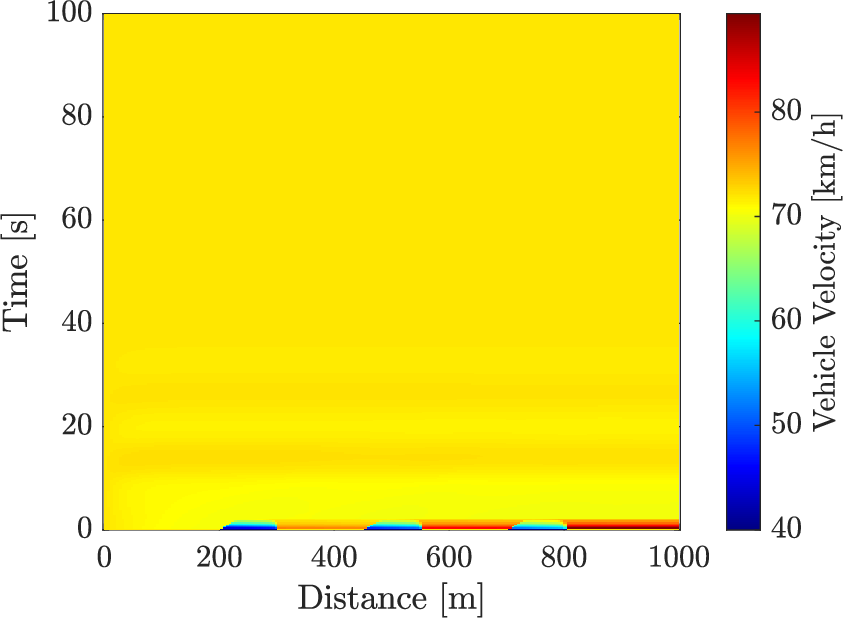}
            \caption{}
            \label{fig:FinCircularRdConVelt_f50S05R}
    \end{subfigure}
    \begin{subfigure}{0.32\textwidth}
            \centering
            \includegraphics[width=\textwidth]{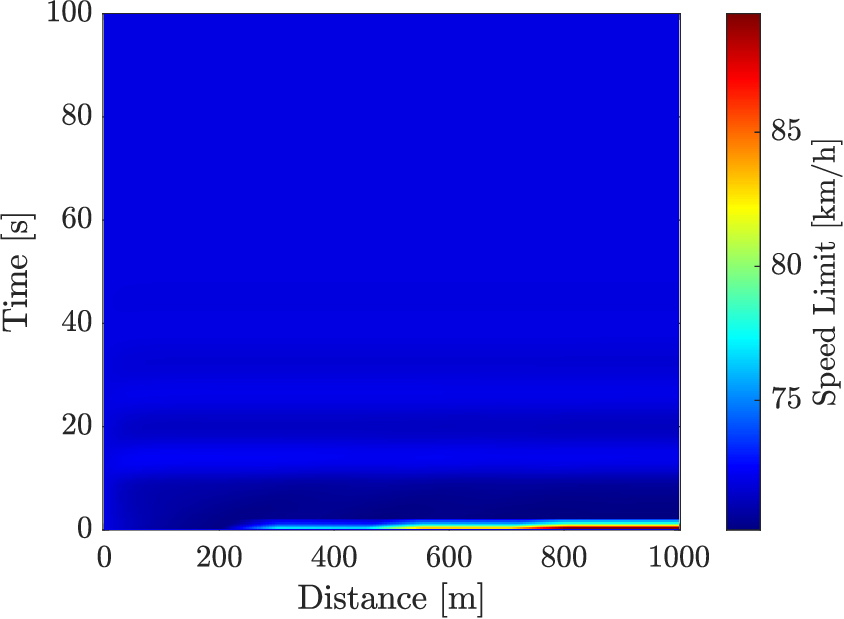}
            \caption{}
            \label{fig:FinCircularRdConLimt_f50S05R}
    \end{subfigure}
    \caption{Density (a), vehicle speed (b), and speed limit (c) for finite-horizon control on a circular road with $t_f = 50$ s, $S = 0.5R$.}
    \label{fig:FinCircularRdt_f50S05R}
\end{figure*}

\cref{fig:FinCircularRdt_f100SR}-\cref{fig:FinCircularRdt_f25SR} show the impact of changing the final time horizon. As the final time decreases, the control action needs to be greater to achieve the desired density faster. This results in speed limits that can become infeasible. Comparing \cref{fig:FinCircularRdt_f50SR}-\cref{fig:FinCircularRdt_f50S05R} it is seen that increasing the value of $S$ from $R$ to $2R$ has minimal affect on the system, but changing it from $R$ to $0.5R$ does, lowering the maximum variable speed limit from over 90 km/h.

\subsubsection{Straight Road Case Study}

\begin{figure*}[ht!]
    \centering
    \begin{subfigure}{0.32\textwidth}
            \centering
            \includegraphics[width=\textwidth]{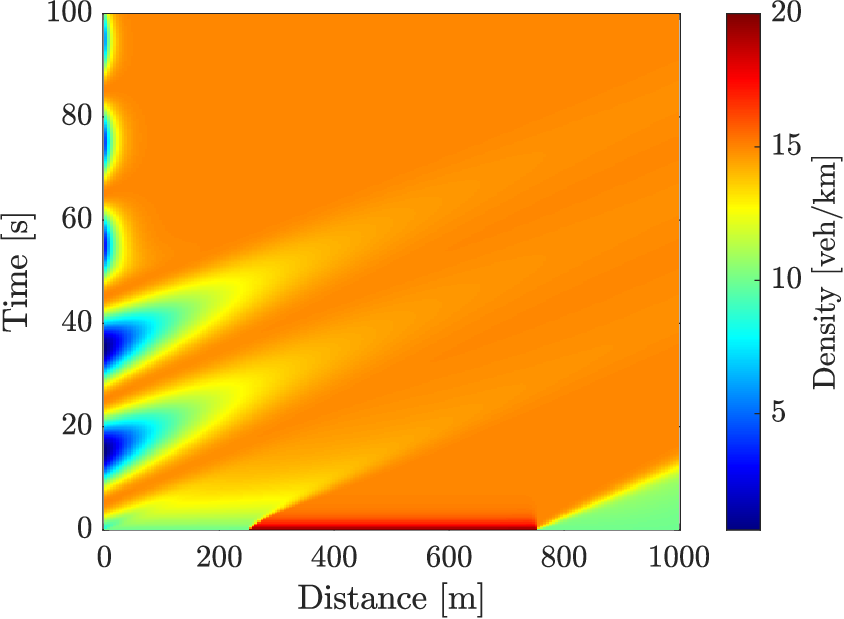}
            \caption{}
            \label{fig:FinStraightRdConRhot_f100SR}
    \end{subfigure}
    \begin{subfigure}{0.32\textwidth}
            \centering
            \includegraphics[width=\textwidth]{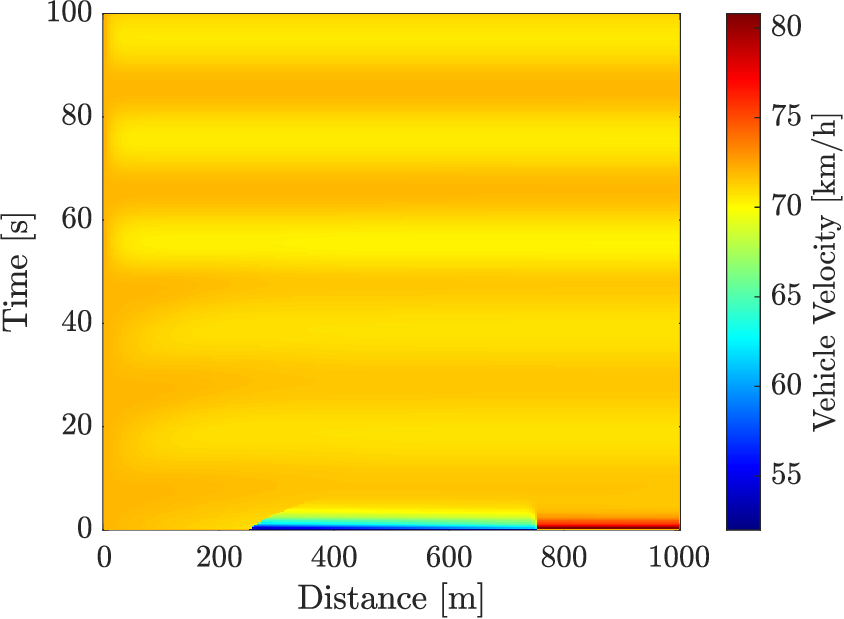}
            \caption{}
            \label{fig:FinStraightRdConVelt_f100SR}
    \end{subfigure}
    \begin{subfigure}{0.32\textwidth}
            \centering
            \includegraphics[width=\textwidth]{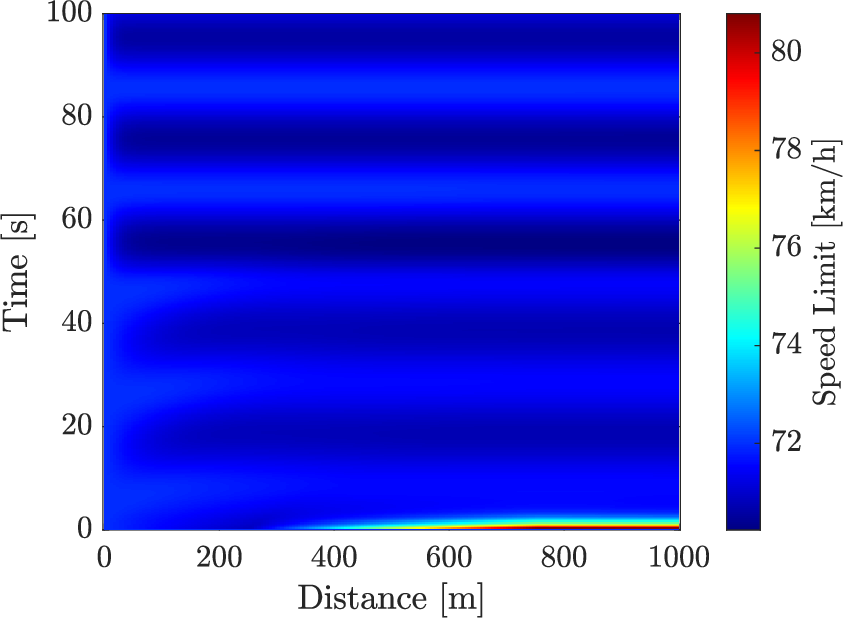}
            \caption{}
            \label{fig:FinStraightRdConLimt_f100SR}
    \end{subfigure}
    \caption{Density (a), vehicle speed (b), and speed limit (c) for finite-horizon control on a straight road with $t_f = 100$ s, $S = R$.}
    \label{fig:FinStraightRdt_f100SR}
\end{figure*}

\begin{figure*}[ht!]
    \centering
    \begin{subfigure}{0.32\textwidth}
            \centering
            \includegraphics[width=\textwidth]{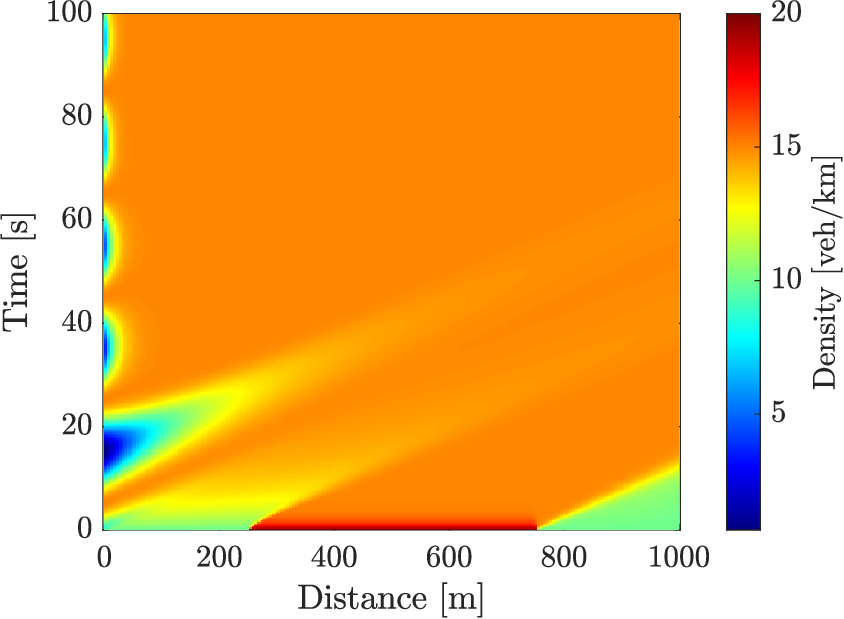}
            \caption{}
            \label{fig:FinStraightRdConRhot_f75SR}
    \end{subfigure}
    \begin{subfigure}{0.32\textwidth}
            \centering
            \includegraphics[width=\textwidth]{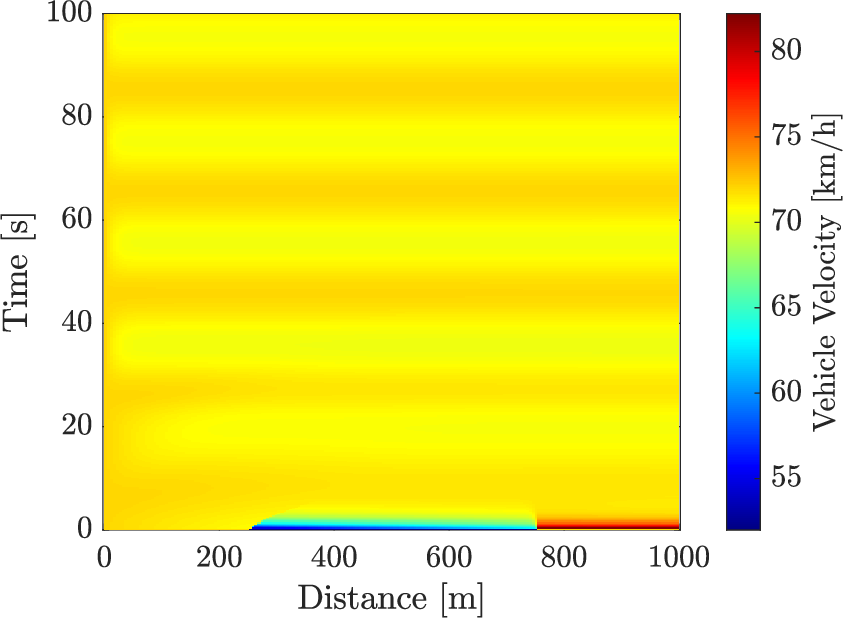}
            \caption{}
            \label{fig:FinStraightRdConVelt_f75SR}
    \end{subfigure}
    \begin{subfigure}{0.32\textwidth}
            \centering
            \includegraphics[width=\textwidth]{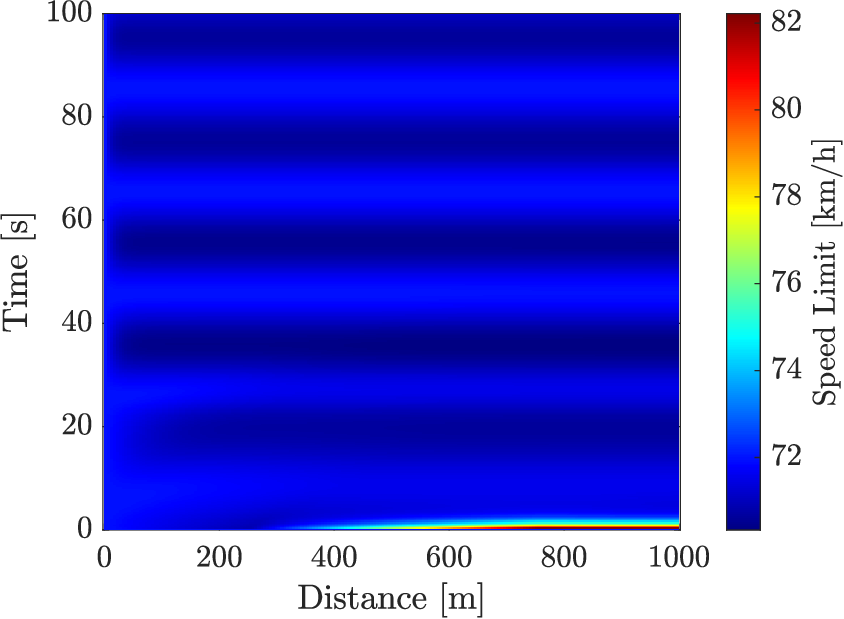}
            \caption{}
            \label{fig:FinStraightRdConLimt_f75SR}
    \end{subfigure}
    \caption{Density (a), vehicle speed (b), and speed limit (c) for finite-horizon control on a straight road with $t_f = 75$ s, $S = R$.}
    \label{fig:FinStraightRdt_f75SR}
\end{figure*}

\begin{figure*}[ht!]
    \centering
    \begin{subfigure}{0.32\textwidth}
            \centering
            \includegraphics[width=\textwidth]{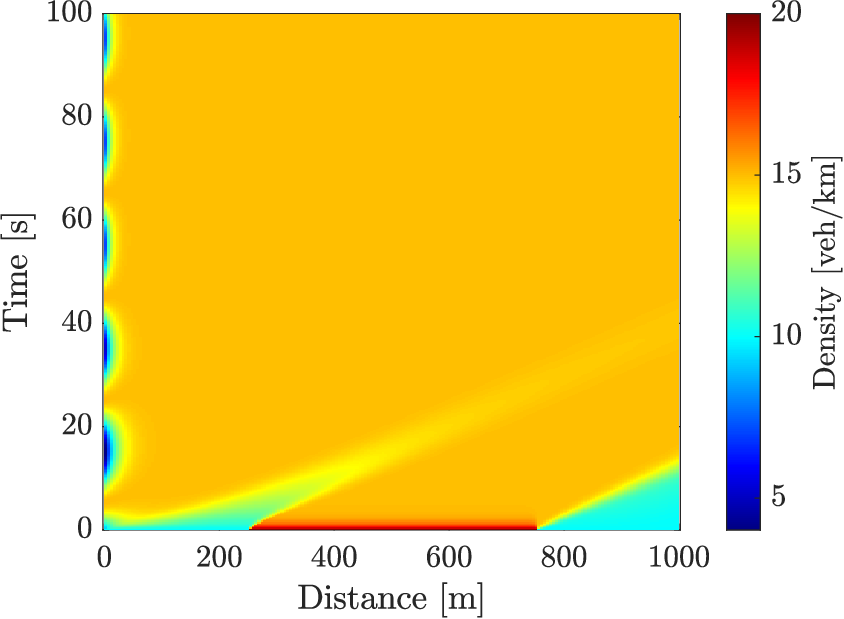}
            \caption{}
            \label{fig:FinStraightRdConRhot_f50SR}
    \end{subfigure}
    \begin{subfigure}{0.32\textwidth}
            \centering
            \includegraphics[width=\textwidth]{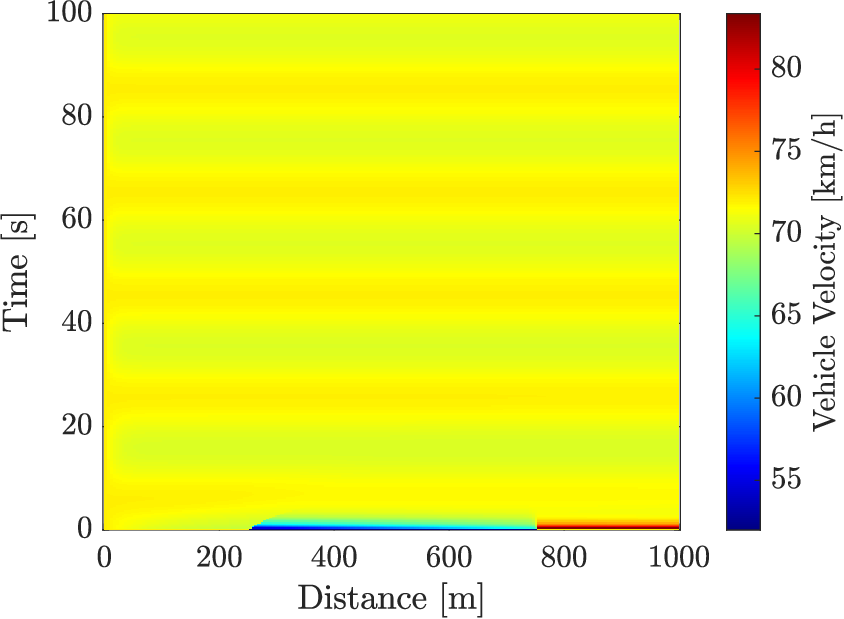}
            \caption{}
            \label{fig:FinStraightRdConVelt_f50SR}
    \end{subfigure}
    \begin{subfigure}{0.32\textwidth}
            \centering
            \includegraphics[width=\textwidth]{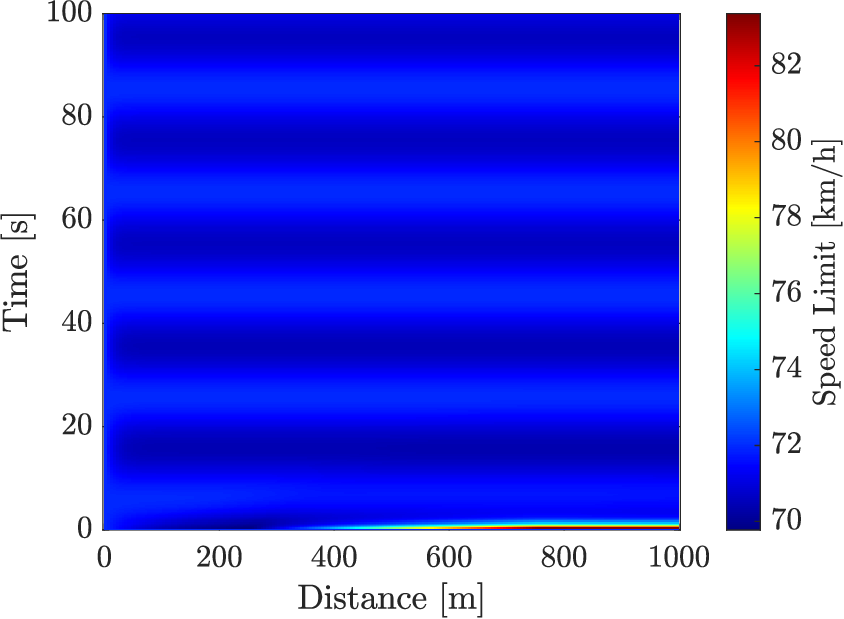}
            \caption{}
            \label{fig:FinStraightRdConLimt_f50SR}
    \end{subfigure}
    \caption{Density (a), vehicle speed (b), and speed limit (c) for finite-horizon control on a straight road with $t_f = 50$ s, $S = R$.}
    \label{fig:FinStraightRdt_f50SR}
\end{figure*}

\begin{figure*}[ht!]
    \centering
    \begin{subfigure}{0.32\textwidth}
            \centering
            \includegraphics[width=\textwidth]{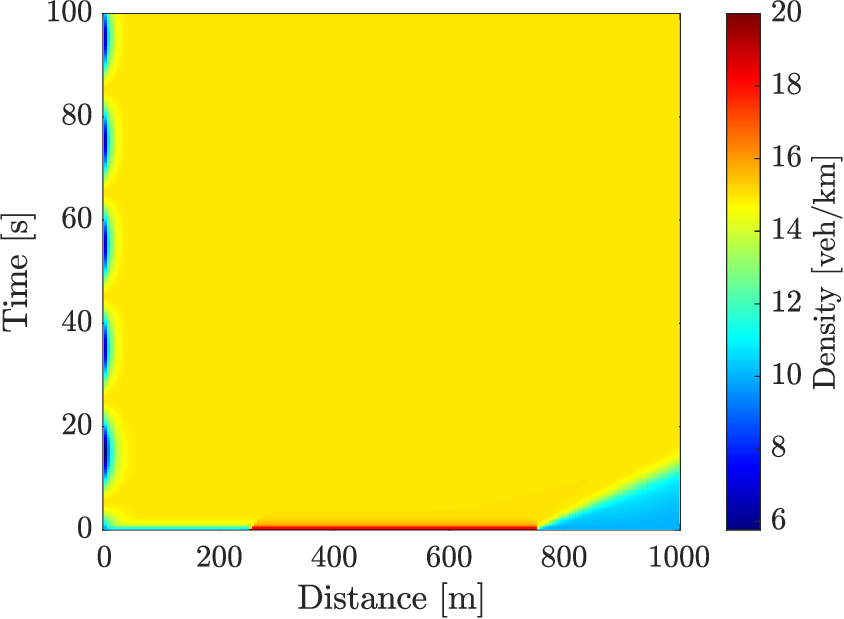}
            \caption{}
            \label{fig:FinStraightRdConRhot_f25SR}
    \end{subfigure}
    \begin{subfigure}{0.32\textwidth}
            \centering
            \includegraphics[width=\textwidth]{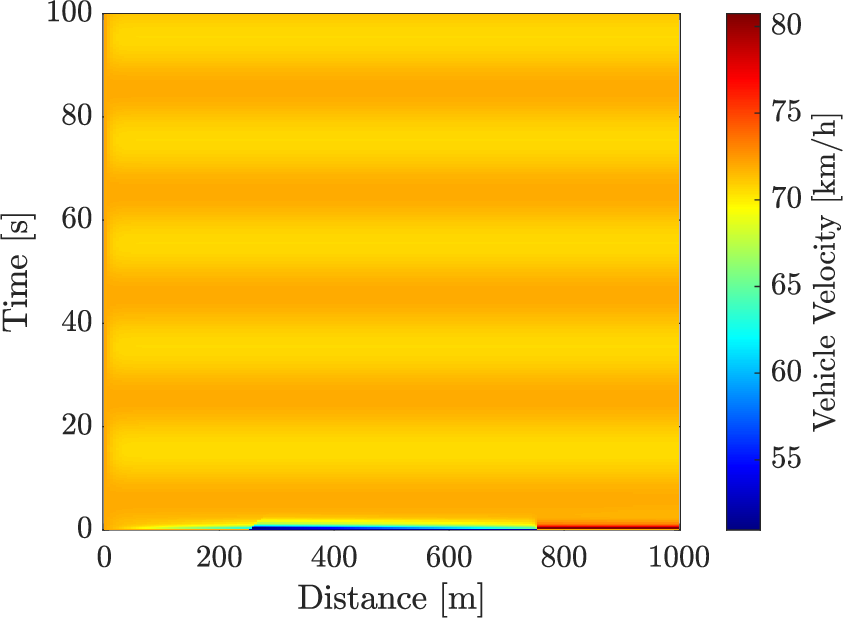}
            \caption{}
            \label{fig:FinStraightRdConVelt_f25SR}
    \end{subfigure}
    \begin{subfigure}{0.32\textwidth}
            \centering
            \includegraphics[width=\textwidth]{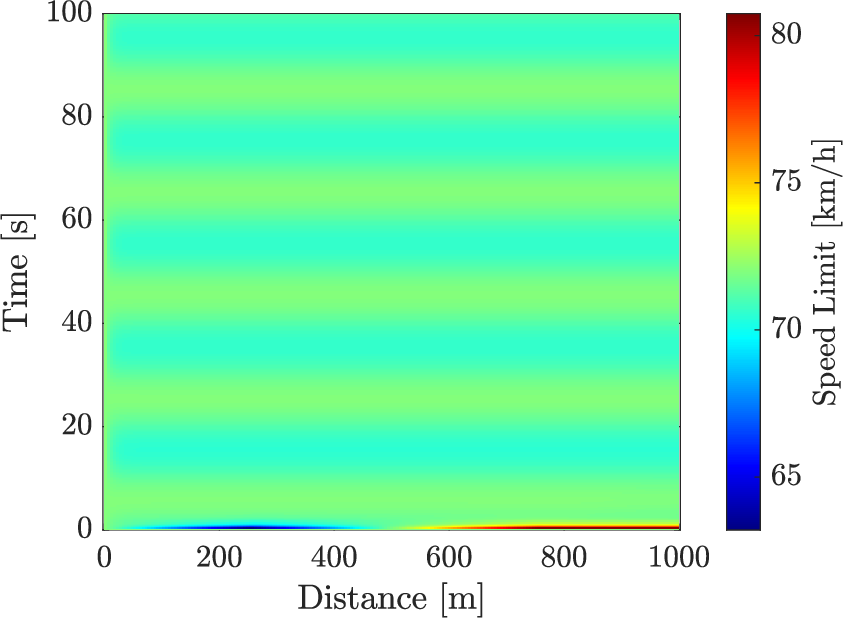}
            \caption{}
            \label{fig:FinStraightRdConLimt_f25SR}
    \end{subfigure}
    \caption{Density (a), vehicle speed (b), and speed limit (c) for finite horizon-control on a straight road with $t_f = 25$ s, $S = R$.}
    \label{fig:FinStraightRdt_f25SR}
\end{figure*}

\begin{figure*}[ht!]
    \centering
    \begin{subfigure}{0.32\textwidth}
            \centering
            \includegraphics[width=\textwidth]{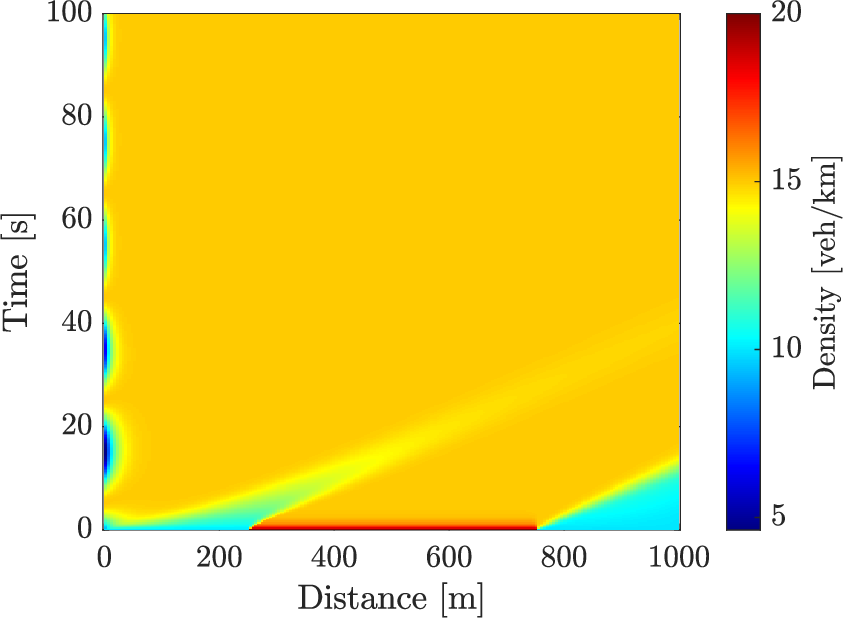}
            \caption{}
            \label{fig:FinStraightRdConRhot_f50S2R}
    \end{subfigure}
    \begin{subfigure}{0.32\textwidth}
            \centering
            \includegraphics[width=\textwidth]{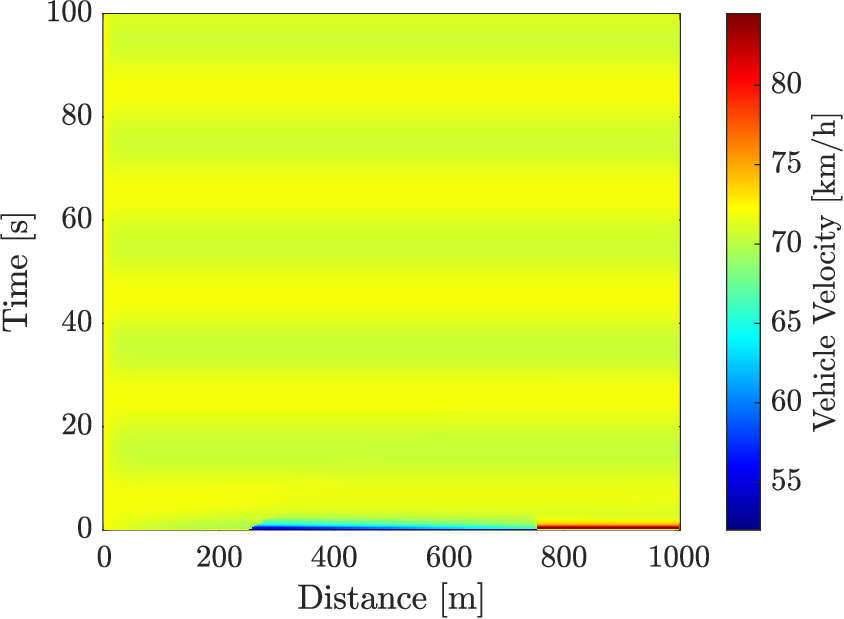}
            \caption{}
            \label{fig:FinStraightRdConVelt_f50S2R}
    \end{subfigure}
    \begin{subfigure}{0.32\textwidth}
            \centering
            \includegraphics[width=\textwidth]{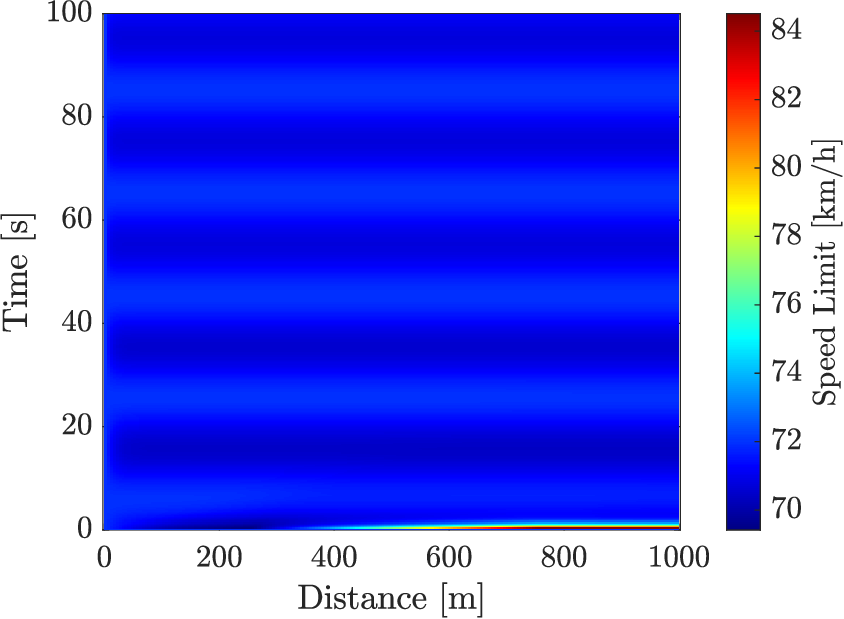}
            \caption{}
            \label{fig:FinStraightRdConLimt_f50S2R}
    \end{subfigure}
    \caption{Density (a), vehicle speed (b), and speed limit (c) for finite-horizon control on a straight road with $t_f = 50$ s, $S = 2R$.}
    \label{fig:FinStraightRdt_f50S2R}
\end{figure*}

\begin{figure*}[ht!]
    \centering
    \begin{subfigure}{0.32\textwidth}
            \centering
            \includegraphics[width=\textwidth]{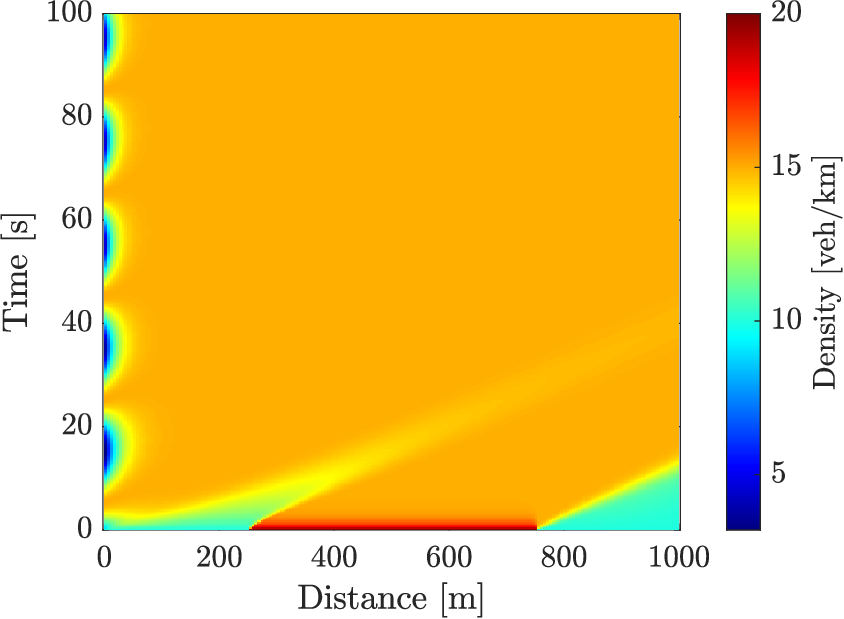}
            \caption{}
            \label{fig:FinStraightRdConRhot_f50S05R}
    \end{subfigure}
    \begin{subfigure}{0.32\textwidth}
            \centering
            \includegraphics[width=\textwidth]{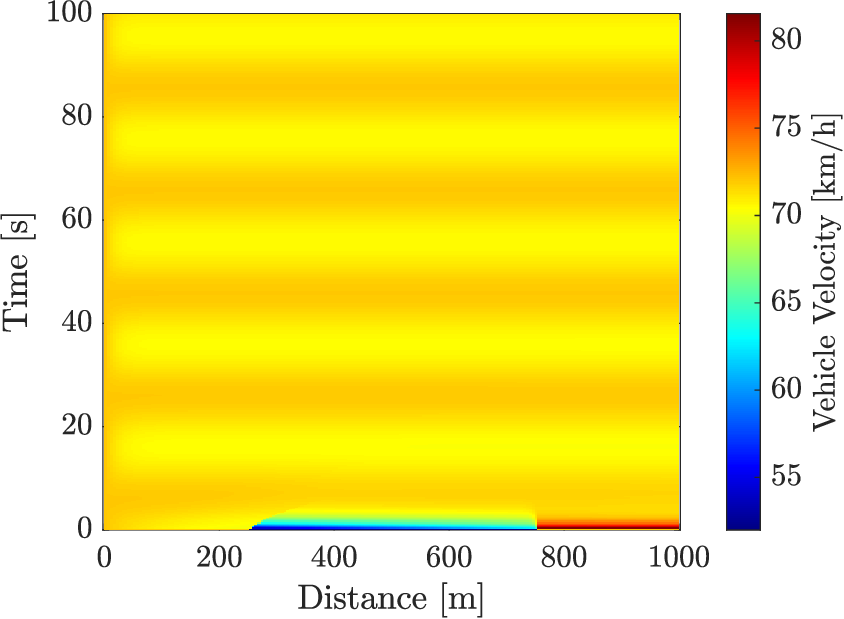}
            \caption{}
            \label{fig:FinStraightRdConVelt_f50S05R}
    \end{subfigure}
    \begin{subfigure}{0.32\textwidth}
            \centering
            \includegraphics[width=\textwidth]{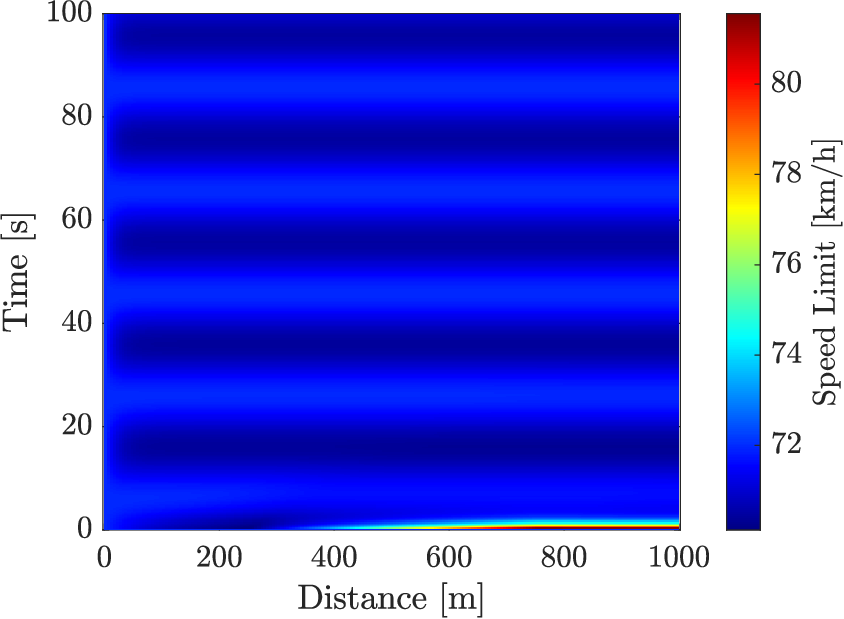}
            \caption{}
            \label{fig:FinStraightRdConLimt_f50S05R}
    \end{subfigure}
    \caption{Density (a), vehicle speed (b), and speed limit (c) for finite-horizon control on a straight road with $t_f = 50$ s, $S = 0.5R$.}
    \label{fig:FinStraightRdt_f50S05R}
\end{figure*}

\cref{fig:FinStraightRdt_f100SR}-\cref{fig:FinStraightRdt_f25SR} show the impact of changing the final time horizon on the straight road scenario with a changing boundary inflow. As expected, as the final time decreases, the variable speed limits can regulate to the desired density faster. In all cases where $t_f<100$ s, the controller does a much better job than the infinite horizon controller at reducing the impact the incoming traffic flow has. Here, changing the value of $S$ does not impact the solution that much as can be seen in \cref{fig:FinStraightRdt_f50S2R} and \cref{fig:FinStraightRdt_f50S05R}.

\section{CONCLUSIONS}\label{sec:conclusion}
This article presented a novel finite horizon LQR approach for the optimal control of traffic using variable speed limits. The LWR traffic model was first linearized and variable speed limit control was integrated by modifying the triangular fundamental diagram. Building on prior work, an infinite horizon state feedback function was extended to accommodate both free flow and congested regimes. A finite horizon LQR problem was formulated, which required the solution of an operator Riccati equation that resulted in a PDE. The analytical solution of the PDE was obtained through the method of characteristics, which gave the state feedback function for the optimal control. The state feedback function was a function of time, space, and the traffic regime.

Then, a sensitivity study was done using both the infinite horizon and finite horizon optimal controllers to investigate the impact the the LQR parameters $Q$ and $R$ as well as the weight on the final cost $S$ and the final time $t_f$ for the finite-horizon control. Both methods were evaluated on mixed traffic on a straight road with changing inflow boundary conditions and a circular road with periodic boundary conditions. It was shown that the finite time horizon controller was able to achieve the goal of regulating density to the desired density target within the prescribed time horizon. Future work will look at developing a feedforward strategies for tracking  time- and space- dependent trajectories, as well as integrating constraints into the problem formulation.



\bibliographystyle{IEEEtran}
\bibliography{references}
\end{document}